\documentclass[aop,citesort,MSNbibl,dvips]{arximspdf}
\usepackage{graphics}
\doi{10.1214/10-AOP550}
\volume{39}
\issue{1}
\pubyear{2011}
\firstpage{104}
\lastpage{138}

\makeatletter

\newcommand{\erf}{\operatorname{erf}}
\newcommand{\Gue}{\operatorname{GUE}}
\newcommand{\Goe}{\operatorname{GOE}}

\newcommand{\e}{\varepsilon}

\newcommand{\R}{\mathbb{R}}
\newcommand{\Z}{\mathbb{Z}}
\newcommand{\ra}{\rightarrow}

\newtheorem{theorem}{Theorem}[section]
\newtheorem{ptheorem}{Partial Theorem}[section]
\newtheorem{lemma}[theorem]{Lemma}
\newtheorem{proposition}[theorem]{Proposition}
\newproclaim{rem}[theorem]{Remark}

\makeatother

\begin{document}
\begin{frontmatter}

\title{Current fluctuations for TASEP: A proof of the Pr\"{a}hofer--Spohn conjecture}
\runtitle{Current fluctuations for TASEP}

\begin{aug}
\author[A]{\fnms{G\'{e}rard} \snm{Ben Arous}\corref{}\thanksref{t1}\ead[label=e1]{benarous@cims.nyu.edu}} and
\author[A]{\fnms{Ivan} \snm{Corwin}\thanksref{t2}\ead[label=e2]{corwin@cims.nyu.edu}}
\runauthor{G. Ben Arous and I. Corwin}
\affiliation{New York University}
\address[A]{Courant Institute\\
\quad of the Mathematical Sciences\\
New York University\\
251 Mercer Street\\
New York, New York 10012\\
USA\\
\printead{e1}\\
\phantom{E-mail: }\printead*{e2}} 
\end{aug}

\thankstext{t1}{Supported in part by NSF Grant
DMS-08-06180.}
\thankstext{t2}{Supported in part by the NSF Graduate
Research Fellowship Program.}

\received{\smonth{5} \syear{2009}}
\revised{\smonth{3} \syear{2010}}

%
\begin{abstract}
We consider the family of two-sided Bernoulli initial conditions for TASEP which,
as the left and right densities ($\rho_-,\rho_+$) are varied, give rise
to shock waves and rarefaction fans---the two phenomena which are
typical to TASEP. We provide a proof of Conjecture 7.1 of
[\textit{Progr. Probab.} \textbf{51} (2002) 185--204]
which characterizes the order of and scaling functions
for the fluctuations of the height function of two-sided TASEP in terms
of the two densities $\rho_-,\rho_+$ and the speed $y$ around which the
height is observed.

In proving this theorem for TASEP, we also prove a fluctuation theorem
for a class of corner growth processes with external sources, or
equivalently for the last passage time in a directed last passage
percolation model with two-sided boundary conditions: $\rho_-$ and
$1-\rho_+$. We provide a complete characterization of the order of and
the scaling functions for the fluctuations of this model's last passage
time $L(N,M)$ as a function of three parameters: the two
boundary/source rates $\rho_-$ and $1-\rho_+$, and the scaling ratio
$\gamma^2=M/N$. The proof of this theorem draws on the results of
[\textit{Comm. Math. Phys.} \textbf{265} (2006) 1--44] and extensively
on the work of [\textit{Ann. Probab.} \textbf{33} (2005) 1643--1697] on
finite rank perturbations of Wishart ensembles in random matrix theory.
\end{abstract}

%
\begin{keyword}[class=AMS]
\kwd{82C22}
\kwd{60K35}.
\end{keyword}
\begin{keyword}
\kwd{Asymmetric simple exclusion process}
\kwd{interacting particle systems}
\kwd{last passage percolation}.
\end{keyword}

\end{frontmatter}

\section{Introduction and results}\label{sec1}
We study the fluctuations of the height function for the Totally
Asymmetric Simple Exclusion Process (TASEP)---a stochastic process of
great interest due to its wide applicability and mathematical
accessibility. Under hydrodynamic scaling, this height function is the
integrated solution to the deterministic Burgers equation
\cite{TL1999s}. This hydrodynamic limit is sensitive to the initial
conditions of TASEP. It is of great interest to determine how the
initial conditions of TASEP affect the random fluctuations of the
height function. Ultimately, one would like to have a dictionary
between initial conditions of TASEP and the resulting orders of the
fluctuations of the height function, along with the scaling functions
and correlation structures. This paper serves to lay some groundwork
for understanding the phenomena which figure into this dictionary. The
two phenomena which must be considered in TASEP are shocks and
rarefaction fans. We study the simplest family of initial conditions
which give rise to both of these phenomena. These initial conditions
are simply Bernoulli independent at each site $x$, with density $\rho
_-$ for $x\leq0$ and $\rho_+$ for $x>0$. We study the fluctuations of
the height function or equivalently the current for these two-sided
initial conditions. We solve an important conjecture of Pr\"
{a}hofer and Spohn \cite{PS2002c} (see also \cite{PLFS2006s}).
Understanding the fluctuation theory for two-sided TASEP provides the
logical link between the well-developed theory for equilibrium initial
conditions ($\rho_-=\rho_+$) \cite{PAFF1994c,PLFS2006s} and step
initial conditions ($\rho_-=1,\rho_+=0$) \cite{KJ2000s}. Two-sided
%
%
\begin{figure}[b]

\includegraphics{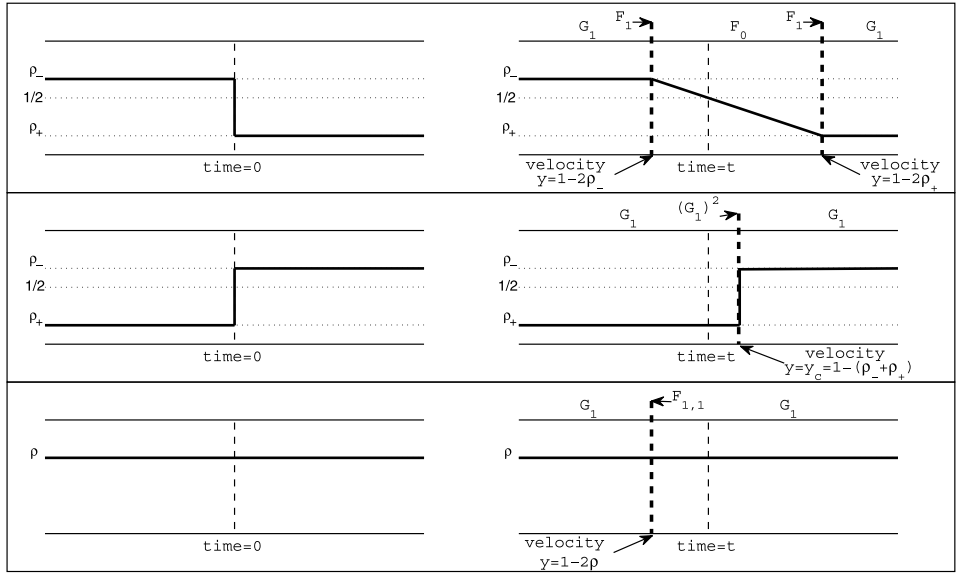}

\caption{Depiction of three types of TASEP (particles move
right) time evolution and identification of different regions of
fluctuations for corresponding height functions. The top diagram
depicts the phenomena of a rarefaction fan. The height function
fluctuations seen by an observer moving at a speed so as to be: outside
of the fan (outside the dashed lines) will be of order $t^{1/2}$ and
Gaussian (denoted $G_1$); inside of the fan (inside the dashed lines)
will be of order $t^{1/3}$ and Tracy--Widom $\Gue$ (denoted $F_0$); on
the edge of the fan (on the dashed lines) will be of order $t^{1/3}$
and Tracy--Widom $\Goe^2$ (denoted $F_1$). Likewise, the middle diagram
depicts the phenomena of a moving shock. The height function
fluctuations seen by an observer moving at a speed so as to be: on the
shock (on the dashed line) will be of order $t^{1/2}$ and ``Gaussian
squared'' [denoted $(G_1)^2$]; away from the shock (off the dashed
line) will be of order $t^{1/2}$ and Gaussian (denoted $G_1$). The
bottom diagram depicts equilibrium initial conditions. The height
function fluctuations seen by an observer moving: at the critical speed
$y_c=1-2\rho$ (on the dashed line) will be of order $t^{1/3}$ (denoted
$F_{1,1}$ and corresponding to what \protect\cite{PLFS2006s} call $F_0$); at
all other speeds (off the dashed line) will be of order $t^{1/2}$ and
Gaussian (denoted $G_1$).}\label{TASEP_fig}
\end{figure}
TASEP interpolates between systems which are in equilibrium and systems
which are entirely out of equilibrium. Our analysis shows how this
interpolation occurs.
The main result, Theorem \ref{thm_conj}, was first conjectured in
\cite{PS2002c} based on a scaling theory and analogous results for the PNG
model and discrete TASEP \cite{BR2000l}. Figure \ref{TASEP_fig}
illustrates the main result of this paper---it shows how the order of
and scaling functions for the fluctuations of the height function for
TASEP depend on the observation location with respect to shocks and
rarefaction fans.

The proof of our main results makes use of the aforementioned result of
\cite{PLFS2006s} for the critical point (equilibrium
$\rho_-=\rho_+=\rho$ and $y=1-2\rho$). For every other set of initial
conditions, the proof relies on the main result of \cite{BBP2005p}, a
paper about the largest eigenvalue of finite rank perturbations of
Wishart (sample covariance) random matrix ensembles. The connection
between these two, seemingly disparate mathematical construction (TASEP
and Wishart ensembles) is due to \cite{KJ2000s} and is facilitated
through an intermediate random process known as directed last passage
percolation (LPP). Our results about fluctuations of currents (or
height functions) for TASEP follow from equivalent results for
fluctuations of the last passage time for a directed last passage
percolation model with two-sided boundary conditions (given in Theorem
\ref{LPP_fluctuations}).

Returning to the model, TASEP is a Markov process $\eta_t$ with state
space $\eta_t \in\{0,1\}^{\Z}$. For a given $t\in\R^{+}$ (time) and
$x\in\Z$ (site location), we say that site $x$ is occupied at time $t$
if $\eta_t(x)=1$ and it is empty if $\eta_t(x)=0$. Given an initial
configuration $\eta_0$ of particles,
the TASEP evolves in continuous time as follows: each particle waits
independent exponentially distributed times and then attempts to
jump one site to its right; if there already exists another
particle in the destination site, the particle does not move
and its waiting time resets (see \cite{TL1999s,TL2005i}
for a rigorous construction of this process). In equilibrium (or
stationary) initial conditions (parametrized by a number $\rho\in
[0,1]$), the $\eta_0(x)$ are independent Bernoulli random variables
with $P(\eta_0(x)=1)=\rho$. In step initial conditions, $\eta_0(x) =1$
for all $x\leq0$ and zero otherwise. Finally, in two-sided initial
conditions (parametrized by a left density $\rho_-$ and a right density
$\rho_+$) $\eta_0(x)$ are independent Bernoulli random variables with
$P(\eta_0(x)=1)=\rho_-$ for $x\leq0$ and $P(\eta_0(x)=1)=\rho_+$
for $x>0$.

A natural and important quantity to study in TASEP is the current of
particles past an observer moving with speed $y$. It is defined as
$J_{yt,t}=$ number of particles to the left of the origin at time
zero and to the right of $yt$ at time $t$ minus number of particles
to the right of the origin at time zero and to the left of $yt$ at
time $t$. The current encodes the same information as the height
function $h_t(j)$ [which we will define in
(\ref{height_function})]:
%
%
\begin{equation} \label{height_to_current}
J_{j,t}=\frac{h_t(j)-j}{2}.
\end{equation}

For equilibrium TASEP with density $\rho$, the law of large numbers and
central limit theorem \cite{PAFF1994c} states that
%
%
\begin{eqnarray}
&&\lim_{t\ra\infty} \frac{J_{yt,t}}{t}=\rho(1-\rho)-y\rho
\qquad
\mbox{almost surely},
\\
%
&&\lim_{t\ra\infty} \frac{J_{yt,t}-E(J_{yt,t})}{\sqrt{t}} = N(0,D_J),
\end{eqnarray}
where $N(0,D_J)$ is a normal with variance
%
%
\begin{equation}
D_J = \rho(1-\rho)|(1-2\rho)-y|.
\end{equation}

For every velocity aside from $y=1-2\rho$, current (and height
function) fluctuations are Gaussian of order $t^{1/2}$. However, for a
single critical velocity the central limit theorem of \cite{PAFF1994c}
is degenerate as the fluctuations are of a lower order than $t^{1/2}$.
In terms of the hydrodynamic limit, this velocity corresponds to the
slope of the characteristic line for Burgers equation. Heuristically
this is the speed at which the initial condition fluctuations travel.
Therefore, at any other speed, the current will depend on more initial
conditions than just that localized to the origin---it is this that
ensures the $t^{1/2}$ fluctuations and Gaussian scaling function for
other velocities. At the critical speed, the initial environment's
fluctuations are of lesser order, and only the dynamic fluctuations
(those due to the actual TASEP process) are felt. These dynamic
fluctuations are of central importance\vspace*{1pt} to understanding KPZ
universality. At the critical speed $y=1-2\rho$ the fluctuations are of
order $t^{1/3}$ and converge, under suitable centering and scaling to a
distribution function
related to the Tracy--Widom $\Gue$ distribution
\cite{PS2002c,PLFS2006s}. Rewriting expression (1.14) of \cite
{PLFS2006s} in terms
of the current $J_{yt,t}$, with $y=1-2\rho$, their $w=0$, and $\chi
=\rho
(1-\rho)$, the result shows that
%
%
\begin{equation}
\lim_{t\ra\infty} P \biggl(\frac{J_{yt,t} -\rho^2
t}{\chi^{2/3}t^{1/3}}\leq x \biggr) = F_{1,1}(x;0;0),
\end{equation}
where $F_{1,1}(x;0;0) = \frac{\partial}{\partial x}
(F_{0}(s)g(x,0))$. The $g(x,0)$ is a scaling function given in their
equation (1.18). See Section \ref{notation} for an overview of
how our notation translates into the notation used in
\cite{PS2002c,PLFS2006s}.

In TASEP starting with step initial conditions, there are no
fluctuations in the initial environment and consequently for every
velocity $y\in(-1,1)$ the current has fluctuations from the dynamics of
order $t^{1/3}$ and with scaling function which corresponds to the
Tracy--Widom $\Gue$ distribution \cite{KJ2000s} (which we write as
$F_0$ so as to be in line with the notation of \cite{BBP2005p}). In
terms of the hydrodynamic limit, the range of speeds $y\in(-1,1)$
corresponds to the entire rarefaction fan, and the fluctuations are
entirely due to the dynamics of TASEP. Ranges of speed bounded away
from the fan correspond to regions which are, in the allotted time,
unchanged by the dynamics of TASEP.

Drawing on the heuristics about the fluctuations along flat and fanned
regions in the hydrodynamic limit, as well as based on a scaling theory
and previous work of \cite{BR2000l} for the PNG model, Pr\"{a}hofer
and Spohn \cite{PS2002c} conjectured that these two fluctuation
theorems (for equilibrium and step initial conditions) arise as cases
of a complete fluctuation theory for two-sided TASEP (see Figure \ref
{TASEP_fig}). In their Conjecture~7.1, Pr\"{a}hofer and Spohn
claimed that the critical point in \cite{PAFF1994c} of $t^{1/3}$
fluctuations for equilibrium TASEP becomes a critical window
(representing the region of the rarefaction fan) as $\rho_-$ is
increased and $\rho_+$ decreased. Ultimately, as $\rho_-=1$ and $\rho
_+=0$ the critical window of velocities
equals the interval $(-1,1)$ as showed in \cite{KJ2000s}.
Likewise, they conjectured Gaussian behavior outside of this window,
as well as in the case where $\rho_-<\rho_+$.

Previous to this paper, part of the Pr\"{a}hofer--Spohn conjecture had
been proved via random matrix techniques in both the papers of Nagao
and Sasamoto \cite{NS2004a} and Baik, Ben Arous and P\'{e}ch\'{e}
\cite{BBP2005p}. Both papers essentially dealt with the case of \mbox
{$\rho_+=0$}
and any $\rho_-\in[0,1]$. Our results are dependent on coupling
arguments which allow us to bootstrap these boundary cases into every
type of two-sided initial condition except for the critical equilibrium
case (which is dealt with via the result of \cite{PLFS2006s}). The
methods of \cite{TS2002d,FR2002c} prove the part of the conjecture
corresponding to the shock ($\rho_-<\rho_+$) by means of a microscopic
Hopf--Lax--Oleinik formula. In these papers, the entire one time
fluctuation process is characterized in the case of the shock. The
scaling conjectured for the rarefaction fan was proved in
\cite{BCS2006c} (in terms of the corresponding corner growth/LPP model
discussed below), though the scaling functions were not addressed therein.

Beyond giving a complete proof of the Pr\"{a}hofer--Spohn conjecture, we
believe that our coupling methods are very natural and provide a highly
intuitive explanation for the transition between Gaussian and
Tracy--Widom scalings. These methods are also useful in studying last
passage percolation models with more general weights and more general
boundary conditions.
In proving the Pr\"{a}hofer--Spohn conjecture,
one may alternatively follow the approach of \cite{PLFS2006s} which is
necessary in the critical case $\rho_-=\rho_+=\rho$ and $y=1-2\rho$.
That argument is very strong and widely applicable. It is based on the
idea of the Schur measure and involves a shift argument and a necessary
analytic continuation argument. Coupling completely avoids these
technical issues and replaces the complex analysis and asymptotic
analysis with simple and intuitive probability. It also seems to be
applicable in certain cases where the Schur measure argument cannot be applied.

Much effort has been devoted to understanding the analogous picture for
ASEP, where particles may move to either the left of the right but are
still subject to the exclusion rule. Progress in this direction was
made in \cite{PAFKS1991m,PAF1992s,PAFF1994c,PAFF1994s,PAFK1995s}
in the early 1990s. The work of Baik and Rains
\cite{BR2000l}, Pr\"{a}hofer and Spohn \cite{PS2002s} and Imamura and
Sasamoto \cite{IS2004f} in the context of the closely related PNG
model and last passage percolation with geometric weights was very
important in formulating and understanding the theory of fluctuations.
Very recently, due to the efforts of Tracy and Widom
\cite{TW2008i,TW2008f,TW2008a,TW2009t,TW2009o}), Derrida and Gerschenfeld
\cite{DG2009c}, Bal\'{a}zs and Sepp\"{a}l\"{a}inen
\cite{BS2008o,BS2008f}, Quastel and Valk\'{o} \cite{QV2007t},
Mountford and Guiol
\cite{MG2005t} significant progress has been made in answering this
question in the general ASEP. In particular, in~\cite{TW2009o}, Tracy
and Widom extend their step initial condition integrable system
approach to ASEP with one-sided Bernoulli initial conditions. In that
case, they observe the exact same fluctuation regimes as for TASEP. At
present, the Pr\"{a}hofer--Spohn conjecture has not been proved for
ASEP. It is tempting to try to use coupling methods to extend the
one-sided picture for ASEP to the two-sided initial condition case. It
is not clear if this is possible, as ASEP is not related to a last
passage percolation model and the coupling occurs at the level of such
a model.

\subsection{Results}

The main result of this paper is a complete proof of \cite{PS2002c}
Conjecture 7.1---our Theorem \ref{thm_conj}.

Following \cite{PS2002c}, assign to a TASEP configuration
$\eta_{t}(j)$ the height function
%
%
\begin{equation}\label{height_function}
h_t(j) = \cases{
\displaystyle2J_{0,t} + \sum_{i=1}^{j} \bigl(1-2\eta_t(i)\bigr),
&\quad$j\geq1$,\cr
\displaystyle2J_{0,t}, &\quad$j=0$,\cr
\displaystyle2J_{0,t} - \sum_{i=j+1}^{0} \bigl(1-2\eta_t(i)\bigr),
&\quad$j\leq-1$.}
\end{equation}
Recall that $J_{0,t}$ is defined as the number of particles
which have crossed the bond $(0,1)$ up to time $t$. For $|y|<1$
denote
%
%
\begin{equation}
\lim_{t\ra\infty} \frac{1}{t} h_t([yt])=\bar{h}(y),
\end{equation}
which exists almost surely due to the law of large numbers
established via the hydrodynamic theory
\cite{HR1981n,TS1998h}. This limit $\bar{h}$ depends not
just on $y$, but also on $\rho_-$ and $\rho_+$ as follows.

If $\rho_-<\rho_+$, then
%
%
\begin{equation}
\bar{h}(y) = \cases{
(1-2\rho_-)y+2\rho_-(1-\rho_-), &\quad for $y\leq y_c$,\cr
(1-2\rho_+)y+2\rho_+(1-\rho_+), &\quad for $y> y_c$,}
\end{equation}
with $y_c = (\rho_+(1-\rho_+)-\rho_-(1-\rho_-))/(\rho_+-\rho_-) =
1-(\rho_-+\rho_+)$.

If $\rho_->\rho_+$, then
%
%
\begin{equation}\quad
\bar{h}(y) =\cases{
(1-2\rho_-)y+2\rho_-(1-\rho_-), &\quad for $y\leq
1-2\rho_-$,\vspace*{2pt}\cr
\frac{1}{2}(y^2+1), &\quad for $1-2\rho_-<y\leq1-2\rho_+$,\vspace
*{2pt}\cr
(1-2\rho_+)y+2\rho_+(1-\rho_+), &\quad for $1-2\rho_+<y$.}
\end{equation}

We prove the following (note that in parenthesis we record the
distribution names as used in \cite{PS2002c}). For an illustration
of the results below, see Figure \ref{TASEP_fig}.
\begin{theorem}[(\cite{PS2002c}, Conjecture 7.1)]\label{thm_conj}

($F_G$). Let either $\rho_-<\rho_+$, $y>y_c$ and
$y<1-\rho_+$, or $\rho_->\rho_+$, $y>1-2\rho_+$ and
$y<1-\rho_+$. Then
%
%
\begin{eqnarray} \label{conj_eq1}
&&\lim_{t\ra\infty} P \bigl(t\bar{h}(y) - h_t([yt])\leq
\bigl(4\rho_+(1-\rho_+)(y-1+2\rho_+)t\bigr)^{1/2} x \bigr) \nonumber\\
[-8pt]\\[-8pt]
&&\qquad= G_1(x).\nonumber
\end{eqnarray}

Let either $\rho_-<\rho_+$, $y<y_c$ and $-\rho_-<y$, or $\rho
_->\rho
_+$, $y<1-2\rho_-$ and \mbox{$-\rho_-<y$}. Then
%
%
\begin{eqnarray}\quad
&&\lim_{t\ra\infty} P \bigl(t\bar{h}(y) - h_t([yt])\leq
\bigl(4\rho_-(1-\rho_-)(-y+1-2\rho_-)t\bigr)^{1/2} x \bigr)\nonumber\\
[-8pt]\\[-8pt]
&&\qquad= G_1(x).\nonumber
\end{eqnarray}

($F_G^2$). Let $\rho_-<\rho_+$ and $y=y_c$, then
%
%
\begin{eqnarray}
&&\lim_{t\ra\infty} P \bigl(t\bar{h}(y)-h_t([yt])\leq
\bigl((\rho_+-\rho_-)t\bigr)^{1/2}x \bigr)\nonumber\\[-8pt]\\[-8pt]
&&\qquad  =
G_1\bigl(\bigl(4\rho_+(1-\rho_+)\bigr)^{-1/2}x\bigr)G_1
\bigl(\bigl(4\rho_-(1-\rho_-)\bigr)^{-1/2}x\bigr).\nonumber
\end{eqnarray}

($F_{\Gue}$). Let $\rho_->\rho_+$ and $1-2\rho_-<y<1-2\rho_+$. Then
%
%
\begin{equation}
\lim_{t\ra\infty} P \bigl(t\bar{h}(y)-h_t([yt])\leq
2^{-1/3}(1-y^2)^{2/3}t^{1/3}x \bigr) =
F_{0}(x).
\end{equation}

($F_{\Goe}^2$). Let $\rho_->\rho_+$ and either $y=1-2\rho_-$ or
$y=1-2\rho_+$. Then
%
%
\begin{equation}\quad
\lim_{t\ra\infty} P \bigl(t\bar{h}(y)-h_t([yt])\leq
2^{-1/3}(1-y^2)^{2/3}t^{1/3}x \bigr) =
F_{1}(x;0).
\end{equation}

($F_0$). Let $\rho_-=\rho=\rho_+$ and $y=1-2\rho$. Then
%
%
\begin{equation}\label{conj_f0}\qquad
\lim_{t\ra\infty} P \bigl(t\bar{h}(y)-h_t([yt])\leq
2^{-1/3}(1-y^2)^{2/3}t^{1/3}x \bigr) =
F_{1,1}(x;0;0).
\end{equation}
\end{theorem}
\begin{rem}\label{rem_conj}
There is one difference between what we prove in Theorem~\ref{thm_conj}
and what is stated in Conjecture 7.1 of
\cite{PS2002c} which is that for the case of the Gaussian
scaling limit, by virtue of the fact that our proof goes by way of a
mapping with directed last passage percolation, there are certain
parts of the Gaussian region (with respect to $y,\rho_-,\rho_+$) for
which our methods do not apply.
\end{rem}

In this model, random weights $w_{i,j}$ are associated to each site
$(i,j)$ in the
upper-right corner of $\Z^{2}$ (with $i\geq0$ and $j\geq0$). The
weights are
usually independent, and often exponential random variables, or
geometric random variables. Every directed (up/right only) path
$\pi$ from $(0,0)$ to $(N,M)$ then has weight $T(\pi)=\sum_{(i,j)\in
\pi} w_{i,j}$, the sum of all weights along the path. The last passage
time from $(0,0)$ to $(N,M)$ is
the maximum path weight over all directed paths:
%
%
\begin{equation}
L(N,M) = \max_{\pi\dvtx(0,0)\ra(N,M)} T(\pi).
\end{equation}
The statistics of $L(N,M)$ are dependent on the choice of
distribution for the random weights, and in certain cases related to
eigenvalue statistics for random matrices.

The current fluctuations for TASEP with step initial conditions were
determined by identifying the height function with a corner growth
model whose growth times correspond with the last passage times for a
specific LPP model with independent rate one exponential
random weights $w_{i,j}$ for $i,j>0$ and boundary weights $w_{i,j}= 0$,
for $i=0$
or $j=0$ \cite{KJ2000s}. Theorem 1.6 of \cite{KJ2000s} shows that as
$M,N\ra\infty$ such that $M/N$ is in a compact subset of $(0,\infty)$,
$L(N,M)$ (as defined by the above weights) is
approximated in distribution by
%
%
\begin{equation}
\bigl(\sqrt{M}+\sqrt{N}\bigr)^2 +
\frac{(\sqrt{M}+\sqrt{N})^{4/3}}{(MN)^{1/6}}\chi_0,
\end{equation}
where $\chi_0$ is distributed as a Tracy--Widom $\Gue$ distribution.
The first term gives the asymptotic average for $L(N,M)$ and the
second term shows that the fluctuations scale like $M^{1/3}$ and
have a well understood scaling function. Via the height functions
mapping, these results translate back into the current fluctuations for
TASEP with step initial conditions. The corner growth model height
function is exactly the random interface bounding the growth region.

This theorem was proved by using the tools of generalized permutations,
the RSK correspondence and Young Tableaux, to relate the distribution
of the last passage time to the distribution of the largest eigenvalue
of a Wishart ensemble, whose
statistics are known to follow the Tracy--Widom $\Gue$ distribution
\cite{KJ2000s}.

By analogy, our method of proof is to first relate two-sided TASEP to
a LPP model, which we appropriately call \textit{LPP with two-sided
boundary conditions} [see (\ref{LPP_weights}) for a
definition], and then to relate the
statistics of the last passage time for that model to the statistics
of eigenvalues of already studied random matrices. The first mapping is already
found in \cite{PS2002c} and relies on Burke's theorem (we review
this mapping in Section \ref{Mapping_section}). LPP with two-sided
boundary conditions is not directly connected to a random matrix
ensemble, however we can realize its last passage time as the
maximum of last passage time for a pair of coupled \textit{LPP
with one-sided boundary conditions} [see (\ref{LPP_weights_one}) for a
definition]. The
last passage time in such one-sided LPP models is related
(see Section 6 of \cite{BBP2005p}) to the largest eigenvalue of
Wishart matrices with
finite rank perturbations. In fact, the phase transitions, with
respect to the magnitude of the finite perturbation, which are
discussed in \cite{BBP2005p} correspond exactly to the transitions
between different orders of and scaling functions for the height
function of two-sided TASEP. By a set of coupling arguments, and
using the results of \cite{BBP2005p} and \cite{PLFS2006s} we provide
a proof of Theorem
\ref{thm_conj}. In proving Theorem \ref{thm_conj}, we reprove the
fluctuation results for step initial conditions as well as for
equilibrium initial conditions (except at the critical point). We also
show that these two results arise from a much more complete picture
(see Figure \ref{TASEP_fig} for an illustration of this).

As noted above, the proof of this theorem relies on understanding the
fluctuations of the last passage time in a LPP model with two-sided
boundary conditions. The specific \textit{LPP with two-sided boundary
conditions} which we will devote much of this paper to studying has
three different types of independent exponential weights $w_{i,j}$:
%
%
\begin{equation}\label{LPP_weights}
w_{i,j} = \cases{
\mbox{exponential of rate } \pi, &\quad if $i>0,j=0$;\cr
\mbox{exponential of rate } \eta, &\quad if $i=0,j>0$;\cr
\mbox{exponential of rate } 1, &\quad if $i>0,j>0$;\cr
\mbox{zero}, &\quad if $i=0,j=0$.}
\end{equation}
In the later part of Section \ref{section_fluctuations}, we will allow
for more general boundary condition where a finite number of columns
and rows can have different (though uniform within the column or row)
rates. The LPP with one-sided boundary conditions is defined similarly
using the weights in (\ref{LPP_weights_one}). At this point,
it is worth remarking that changing the distribution of a finite number
of weights does not have any affect on the asymptotic fluctuations of
the last passage time (see Lemma \ref{finite_weight_changes}).


The statistics considered in this paper are the last passage times
$L_2(N,M)$ (we use a subscript 2 to denote two-sided), from
$(0,0)$ to $(N,M)$, as $N$ and $M$ go to infinity together such that
$M/N=\gamma^2$. Note that $N$ denotes the number of columns and $M$ the
number of rows. Such statistics can be parametrized in terms of the two
boundary condition rates $\pi$ and $\eta$, as well as the scaling
parameter $\gamma$. It is worth keeping in mind that the boundary rates
$\pi$ and $\eta$ correspond with the TASEP densities $\rho_-$ and
$1-\rho_+$, and the scaling parameter $\gamma$ corresponds (in a
slightly more complicated way) with the TASEP velocity $y$.

With this connection in mind, we completely characterize both the order
and the scaling functions for the fluctuations of the last passage time
of LPP with two-sided boundary conditions in terms of the three
parameters $\pi$, $\eta$ and $\gamma$. As noted before, the main
result we appeal to in this paper is from \cite{BBP2005p} (extended to
the case $\gamma<1$ in~\cite{AO2008t}) which classifies the
fluctuations of the largest eigenvalue of complex Wishart ensembles
with finite rank perturbations. There is a single critical point which
does not yield to our method of argument, but this corresponds exactly
with the critical point considered in \cite{PLFS2006s}. Using these
%
%
\begin{figure}

\includegraphics{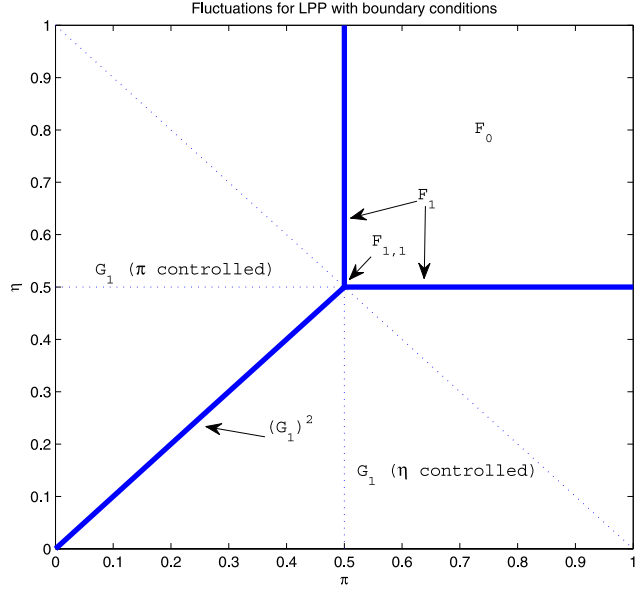}

\caption{Fluctuation diagram for $\gamma=1$. Note that all $G_1$ and
$(G_1)^2$ regions have $M^{1/2}$ order fluctuations while all other
regions have $M^{1/3}$ order fluctuations.}\label{gamma_1}
\end{figure}
%
%
\begin{figure}

\includegraphics{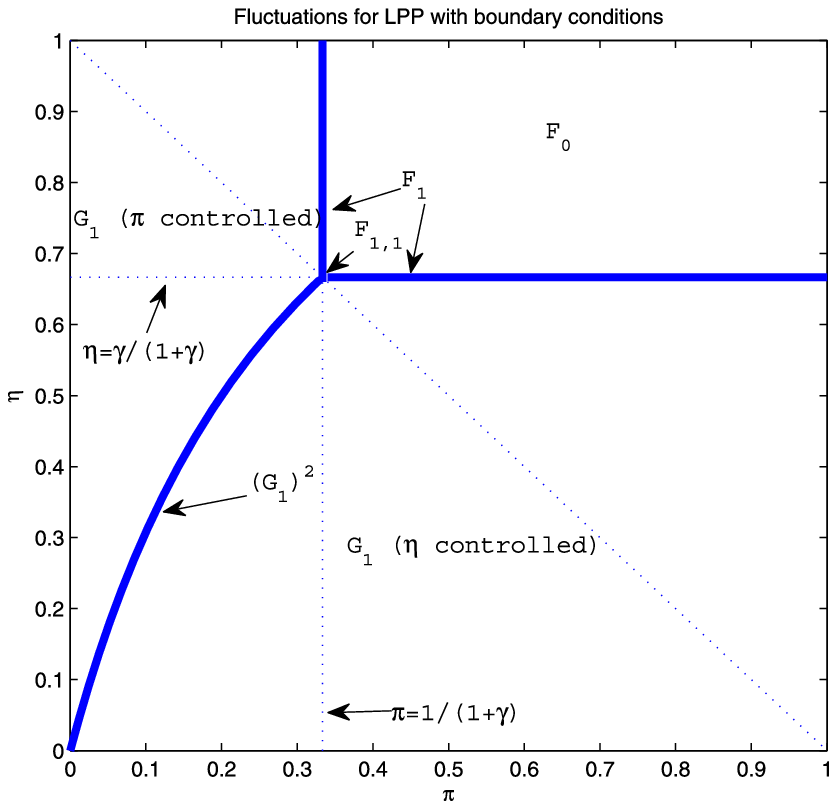}

\caption{Fluctuation diagram for $\gamma=2$. Compared with Figure
\protect\ref{gamma_1}, the effect of changing $\gamma$ is that the
region of
$M^{1/3}$ fluctuations has shifted up and to the left along the
anti-diagonal.}\label{gamma_2}
\end{figure}
two results and coupling arguments, we prove our LPP with two-sided
boundary conditions classification theorem (see Figures \ref{gamma_1}--\ref{gamma_1_over_2}).
\begin{theorem} \label{LPP_fluctuations}
\begin{enumerate}[(2)]
\item[(1)]
For $\gamma\in(0,\infty)$ and $M/N\ra\gamma^2$, then for
$\pi,\eta$ such that $\pi>\frac{1}{1+\gamma}$ and
$\eta>\frac{\gamma}{1+\gamma}$ (the $\Gue$ region)
%
%
\begin{equation}
P\biggl(L_2(N,M)\leq(1+\gamma^{-1})^2M +\frac{(1+\gamma
)^{4/3}}{\gamma}M^{1/3}x\biggr)\ra F_0(x),
\end{equation}
where $F_0(x)$ is the Tracy--Widom $\Gue$ distribution function.
\item[(2)]
For $\gamma\in(0,\infty)$ and $M/N\ra\gamma^2$, then for
$\pi,\eta$ such that $\pi>\frac{1}{1+\gamma}$ and
$\eta=\frac{\gamma}{1+\gamma}$ or $\pi=\frac{1}{1+\gamma}$ and
$\eta>\frac{\gamma}{1+\gamma}$ (the $\Goe^2$ region),
%
%
\begin{equation}
P\biggl(L_2(N,M)\leq(1+\gamma^{-1})^2M
+\frac{(1+\gamma)^{4/3}}{\gamma}M^{1/3}x\biggr)\ra F_1(x),
\end{equation}
where $F_1(x)$ is the square of the Tracy--Widom $\Goe$ distribution
function.
\item[(3)]
For $\gamma\in(0,\infty)$ and $M/N\ra\gamma^2$, then for
$\pi=1/(1+\gamma)$ and $\eta=\gamma/(1+\gamma)$,
%
%
\begin{equation}
P\biggl(L_2(N,M)\leq(1+\gamma^{-1})^2M
+\frac{(1+\gamma)^{4/3}}{\gamma}M^{1/3}x\biggr)\ra F_{1,1}(x;0;0),
\end{equation}
where $F_{1,1}(x;0;0)$ is the same distribution as what
\cite{PLFS2006s} refer to as $F_0(x)$.
\item[(4)]
For $\gamma\in(0,\infty)$ and $M/N\ra\gamma^2$, then for
$\pi,\eta$ such that $\pi<1/(1+\gamma)$ and
$\eta>\frac{\pi}{\pi(1-\gamma^{-2})+\gamma^{-2}}$ [the $G$ ($\pi$
controlled) region],
%
%
\begin{eqnarray}
&&P \biggl(L_2(N,M)\leq
\biggl(\pi^{-1}+\frac{\pi^{-1}\gamma^{2}}{\pi^{-1}-1} \biggr)N +
\biggl(\pi^{-2}-\frac{\pi^{-2}\gamma^{2}}{(\pi^{-1}-1)^2} \biggr)^{1/2}
N^{1/2} x \biggr)\nonumber\\[-8pt]\\[-8pt]
&&\qquad\ra G_1(x),\nonumber
\end{eqnarray}
where $G_1(x)=\erf(x)$.

%
%
\begin{figure}

\includegraphics{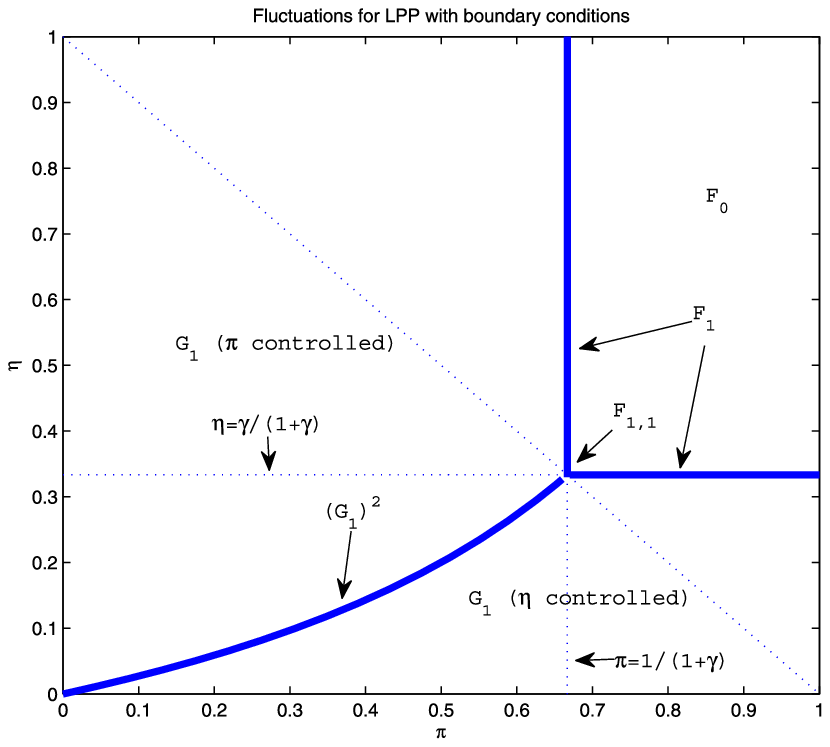}

\caption{Fluctuation diagram for $\gamma=0.5$. Compared with Figure
\protect\ref{gamma_1}, the effect of changing $\gamma$ is that the
region of
$M^{1/3}$ fluctuations has shifted down and to the right along the
anti-diagonal.}\label{gamma_1_over_2}
\end{figure}

Likewise for $\pi,\eta$ such that $\eta<\gamma/(1+\gamma)$ and
$\eta<\frac{\pi}{\pi(1-\gamma^{-2})+\gamma^{-2}}$ [the $G$
($\eta$~controlled) region],
%
%
\begin{eqnarray}
&&
P \biggl(L_2(N,M)\leq
\biggl(\eta^{-1}+\frac{\eta^{-1}\gamma^{-2}}{\eta^{-1}-1} \biggr)M +
\biggl(\eta^{-2}-\frac{\eta^{-2}\gamma^{-2}}{(\eta^{-1}-1)^2} \biggr)^{1/2}
M^{1/2} x \biggr)\nonumber\\[-8pt]\\[-8pt]
&&\qquad\ra G_1(x),\nonumber
\end{eqnarray}
where $G_1(x)=\erf(x)$.
\item[(5)]
For $\gamma\in(0,\infty)$ and $M/N\ra\gamma^2$, then for
$\pi,\eta$ such that $\pi+\eta<1$ and
$\eta=\frac{\pi}{\pi(1-\gamma^{-2})+\gamma^{-2}}$ (the $G^2$ line),
%
%
\begin{eqnarray}
&&
P \biggl( L_2(N,M)\leq\biggl(\pi^{-1}+\frac{\pi^{-1}\gamma^{2}}{\pi
^{-1}-1} \biggr)N\nonumber\\[-8pt]\\[-8pt]
&&\qquad\quad\hspace*{40.8pt}{} + \biggl(\frac{(1-\pi+\pi\gamma^2)((1-\pi)^2-\pi
^2\gamma^2)}{\gamma^2\pi^2(1-\pi)^2} \biggr)^{1/2} N^{1/2} x \biggr)\nonumber\\
&&\qquad\ra G_1\biggl(x\frac{\gamma}{\sqrt{1-\pi+\gamma^2 \pi}}\biggr)G_1\biggl(x
\frac{\sqrt{1-\pi+\gamma^2 \pi}}{\gamma}\biggr).
\end{eqnarray}
\end{enumerate}
\end{theorem}

It is worth noting that there are many ways to write the expressions
above, and our choices are to facilitate the greatest ease in our
proofs.

This type of fluctuation classification picture has been previously
discussed in \cite{BR2000l} and \cite{Sa2007f}. In fact, in
\cite{BR2000l} Baik and Rains provide a proof of an analogous fluctuation
classification result for two closely related particle system models:
LPP with geometric weights, and the polynuclear growth model. Recently,
\cite{BPS2009t} studied two-speed (different though related to
two-sided initial conditions, half flat and half Bernoulli) TASEP and
proved a fluctuation classification theorem for that model. As noted
before, \cite{BCS2006c} previously provided the order of fluctuations
for LPP models with two-sided boundary condition corresponding to the
rarefaction fan. With an even more general type of boundary condition,
the paper establishes $t^{1/3}$ scaling for the fluctuations of the
last passage time.

\subsection{Notation}\label{notation}

In a paper such as this which connects two different lines of
thought, it is easy to become lost in the disparity between
notations. We will adopt notation in the style of \cite{BBP2005p}
throughout, and when making connections with distributions as found
in papers such as \cite{PS2002c,PLFS2006s,BR2000l},
we will take care to make note of the alternative
notation used in those contexts. In this section, we define all of
the distributions which we will encounter herein and provide
references for their previous use and definition.

\begin{enumerate}[(2)]
\item[(1)] $G_k(x)$ is a family of distributions defined in
\cite{BBP2005p}, Definition 1.2 and Lem\-ma~1.1. It represents the distribution of the
largest eigenvalue of a $k\times k$ $\Gue$. From this representation,
it is clear that $G_1(x)=\erf(x)$, the standard Gaussian distribution function.
\item[(2)] $F_J(x;x_1,\ldots, x_J)$ is a family of distributions defined in
\cite{BBP2005p}, Definition 1.3. In the case when the $x_j=0$ for all
$j$, these distributions coincide with those from \cite{BBP2005p},
Definition 1.1. Of note is $F_0(x)$ which is often written as $F_{\Gue
}$, the $\Gue$ Tracy--Widom distribution function, and $F_1(x;0)$ which
is often written as $F_{\Goe}(x)^2$, where $F_{\Goe}$ is the $\Goe$
Tracy--Widom distribution function.
\item[(3)] $F_{J,I}(x;x_1,\ldots, x_J;y_1,\ldots, y_I)$ is a family of
distributions which we conjecture come up in LPP with thick two-sided
boundary conditions. The only member of this family for which we know
the correct definition is $F_{1,1}(x;0;0)$ which corresponds to the
distribution denoted by $F_0$ in \cite{PLFS2006s}. As of yet, we do
not know how the other distributions should be defined.
\end{enumerate}

\subsection{Outline}
The main theorems (Theorems \ref{thm_conj} and \ref
{LPP_fluctuations}) have already been recorded above in this section.
Section \ref{section_fluctuations} provides an intuitive sketch of
the proof for Theorem \ref{LPP_fluctuations}. Section
\ref{Mapping_section} explains the connection between the LPP with
two-sided boundary conditions and the TASEP with two-sided initial
conditions as well as briefly sketches how to translate the result
of Theorem \ref{LPP_fluctuations} into a proof of Theorem
\ref{thm_conj}. Section \ref{proofs} gives the full proof
of the two main theorems, complete with the necessary technical
lemmas for the coupling arguments.

\section{Fluctuations in last passage percolation with boundary
conditions} \label{section_fluctuations}
We start this section by reviewing the result of \cite{BBP2005p}
which relates directed last passage percolation with boundary
conditions to finite rank perturbations of Wishart ensembles. We
then apply these results to prove Theorem \ref{LPP_fluctuations}
which fully characterizes the fluctuations of last passage times in
terms of boundary conditions and the ratio $M/N =\gamma^2$. Using
coupling arguments, supplemented in one case by the result of
\cite{PLFS2006s}, we provide both the order and the scaling function
for these fluctuations. Using the exact same arguments but fully
taking advantage of the scope of the results of \cite{BBP2005p}, we
prove almost all of the cases in Partial Theorem
\ref{LPP_fluctuations_thick}. In this section, we will only sketch
our proofs, which can be found in entirety in Section~\ref{proofs}.

\subsection{LPP with one-sided boundary conditions}

Consider a directed last passage percolation model with one-sided
boundary conditions defined as follows:
%
%
\begin{equation}\label{LPP_weights_one}
w_{i,j} = \cases{
\mbox{exponential of rate } \eta, &\quad if $i=0,j>0$,\cr
\mbox{exponential of rate } 1, &\quad if $i>0,j>0$,\cr
\mbox{zero}, &\quad if $i\geq0,j=0$.}
\end{equation}
Let $L_1(N,M)$ denote the last passage time from $(0,0)$ to $(N,M)$
(for this LPP model with one-sided boundary conditions, but also for
any LPP model with thick one-sided boundary conditions). Then the
distribution of $L_1(N,M)$ is related to the distribution of the
largest eigenvalue of the normalized covariance matrix
$\frac{1}{M}XX'$ where $X$ is $N\times M$ and each column is drawn
(independent of other columns) from a complex $N$-dimensional
Gaussian distribution with covariance matrix $\Sigma$. The matrix
$\Sigma$ has eigenvalues all equal to one aside from a single one,
which is $l_1=\eta^{-1}$. Depending on the value of $\eta^{-1}$,
$L_1(N,M)$ behaves differently.

The following theorem is adapted from Theorem 1.1 of
\cite{BBP2005p} and the extension to all $\gamma\in(0,\infty)$
given in \cite{AO2008t}, as applied to the one-sided boundary
condition LPP. The connection between the largest eigenvalue and the
LPP with one-sided boundary conditions given above is explained in
Section 6 of \cite{BBP2005p} and is briefly rehashed in Remark
\ref{BBP_LPP_rem}.
\begin{proposition}\label{LPP_BBP}
With $L_1(N,M)$ defined as above, as $M,N\ra\infty$ while $M/N=\gamma
^2$ is in a compact subset of $(0,\infty)$, the following hold for any
real $x$ in a compact set.

\begin{enumerate}[(2)]
\item[(1)] When $\eta>\frac{\gamma}{1+\gamma}$,
%
%
\begin{equation}
P \biggl(L_1(N,M)\leq(1+\gamma^{-1})^2 M +
\frac{(1+\gamma)^{4/3}}{\gamma}M^{1/3}x \biggr)\ra F_0(x);
\end{equation}

\item[(2)] When $\eta=\frac{\gamma}{1+\gamma}$,
%
%
\begin{equation}
P \biggl(L_1(N,M)\leq(1+\gamma^{-1})^2 M +
\frac{(1+\gamma)^{4/3}}{\gamma}M^{1/3}x \biggr)\ra F_1(x);
\end{equation}

\item[(3)] When $\eta<\frac{\gamma}{1+\gamma}$,
%
%
\begin{eqnarray}
&&P \biggl(L_1(N,M)\leq
\biggl(\eta^{-1}+\frac{\eta^{-1}\gamma^{-2}}{\eta^{-1}-1} \biggr)M +
\biggl(\eta^{-2}-\frac{\eta^{-2}\gamma^{-2}}{(\eta^{-1}-1)^2} \biggr)^{1/2}
M^{1/2} x \biggr)\nonumber\\[-8pt]\\[-8pt]
&&\qquad\ra G_1(x).\nonumber
\end{eqnarray}
\end{enumerate}
\end{proposition}
\begin{rem}\label{BBP_LPP_rem}
The connection between last passage time in LPP with one-sided
boundary conditions and the largest eigenvalue of the spiked\break Wishart
ensemble was observed in \cite{BBP2005p}. The connection is not via
an exact map but rather an equality of
distributions. Proposition 6.1 of \cite{BBP2005p} records this fact
and explains how a modification of the argument in \cite{KJ2000s}
can be used to prove this.
\end{rem}

An intuitive explanation for the cutoff of $\eta^{-1}=1+\gamma^{-1}$
in terms of a simple calculus problem of maximizing the law of large
numbers for LPP paths forced to travel a specific fraction of the
way along the left column can be found in Section~6 of
\cite{BBP2005p}.

\subsection{LPP with two-sided boundary conditions}

Presently, we turn our attention to the LPP with two-sided boundary
conditions as defined in (\ref{LPP_weights}): on the
left-most column there are exponential weights of rate $\eta$, on
the bottom-most row there are exponential weights of rate $\pi$ at
the origin there is a weight of zero, and for all strictly positive
lattice points the weight is of rate one. Define, respectively,
$X(N,M)$ and $Y(N,M)$ as the coupled last passage times of paths
which have taken the first step to the right and the first step up
(resp.). One should be careful to note that we are not
conditioning on the location of the optimal path, but rather, for
each configuration of weights, defining $X$ to be the length of the
optimal path which first goes right, and $Y$ the length of the
optimal path which first goes up. It is clear then that $X$ and $Y$
are coupled, dependent and that
%
%
\begin{equation}
L_2(N,M) = \max(X(N,M),Y(N,M) ).
\end{equation}
Consider now the marginals of $X$ and $Y$ and observe that each of
these marginals is of the type of the last passage time for a LPP model
with one-sided boundary conditions. The boundary conditions for $Y$ are
exactly as above ($\eta$ weights and an $N$ by $M$ region). However,
for $X$, the boundary conditions are $\pi$ weights and an $M$ by $N$
region (note that the region has been flipped in order to conform with
the setup for Proposition \ref{LPP_BBP}). From this observation, we can
apply Proposition \ref{LPP_BBP} to completely characterize the
marginals of the joint distribution for the pair $(X,Y)$. Note that
while $X$ and $Y$ are not exactly of the form of a last passage time
for a LPP with one-sided boundary conditions, they only differ by a
finite number of weights and therefore have the exact same asymptotic
statistics via Lemma \ref{finite_weight_changes}.
\begin{proposition}\label{X_dist}
With $X(N,M)$ defined as above, as $M,N\ra\infty$ while $M/N=\gamma
^2$ is in a compact subset of $(0,\infty)$, the following hold for any
real $x$ in a compact set:

\begin{enumerate}[(2)]
\item[(1)] when $\pi>\frac{1}{1+\gamma}$,
%
%
\begin{equation}
P \biggl(X(N,M)\leq(1+\gamma)^2 N +
\frac{(1+\gamma^{-1})^{4/3}}{\gamma^{-1}}N^{1/3}x \biggr)\ra F_0(x);
\end{equation}

\item[(2)] when $\pi=\frac{1}{1+\gamma}$,
%
%
\begin{equation}
P \biggl(X(N,M)\leq(1+\gamma)^2 N +
\frac{(1+\gamma^{-1})^{4/3}}{\gamma^{-1}}N^{1/3}x \biggr)\ra F_1(x);
\end{equation}

\item[(3)] when $\pi<\frac{1}{1+\gamma}$,
%
%
\begin{eqnarray}
&&
P \biggl(X(N,M)\leq
\biggl(\pi^{-1}+\frac{\pi^{-1}\gamma^{2}}{\pi^{-1}-1} \biggr)N +
\biggl(\pi^{-2}-\frac{\pi^{-2}\gamma^{2}}{(\pi^{-1}-1)^2} \biggr)^{1/2}
N^{1/2} x \biggr)\nonumber\\[-8pt]\\[-8pt]
&&\qquad\ra G_1(x).\nonumber
\end{eqnarray}
\end{enumerate}
\end{proposition}
\begin{proposition}\label{Y_dist}
With $Y(N,M)$ defined as above, as $M,N\ra\infty$ while
$M/N=\gamma^2$ is in a compact subset of $(0,\infty)$, the following
hold for any real $x$ in a compact set:

\begin{enumerate}[(2)]
\item[(1)] when $\eta>\frac{\gamma}{1+\gamma}$,
%
%
\begin{equation}
P \biggl(Y(N,M)\leq(1+\gamma^{-1})^2 M +
\frac{(1+\gamma)^{4/3}}{\gamma}M^{1/3}x \biggr)\ra F_0(x);
\end{equation}

\item[(2)] when $\eta=\frac{\gamma}{1+\gamma}$,
%
%
\begin{equation}
P \biggl(Y(N,M)\leq(1+\gamma^{-1})^2 M +
\frac{(1+\gamma)^{4/3}}{\gamma}M^{1/3}x \biggr)\ra F_1(x);
\end{equation}

\item[(3)] when $\eta<\frac{\gamma}{1+\gamma}$,
%
%
\begin{eqnarray}
&&P \biggl(Y(N,M)\leq
\biggl(\eta^{-1}+\frac{\eta^{-1}\gamma^{-2}}{\eta^{-1}-1} \biggr)M +
\biggl(\eta^{-2}-\frac{\eta^{-2}\gamma^{-2}}{(\eta^{-1}-1)^2} \biggr)^{1/2}
M^{1/2} x \biggr)\nonumber\\[-8pt]\\[-8pt]
&&\qquad\ra G_1(x).\nonumber
\end{eqnarray}
\end{enumerate}
\end{proposition}

We assume that both $\eta$ and $\pi$ are between zero and one.
In fact, it is clear from our proofs that the order and fluctuation
of our two-sided last passage time $L_2(N,M)$ for parameters $\eta$
and $\pi$ is that same as that for parameters $\eta\wedge
1,\pi\wedge1$. Thus, it suffices to consider only $\eta,\pi\in
[0,1]^2$.

For each of $X$ and $Y$, there are two regions of different
fluctuation orders, and one critical point which has $1/3$ order
fluctuations. We call the point
$(\pi,\eta)=(1/(1+\gamma),\gamma/(1+\gamma))$ the critical point for
the pair $\pi,\eta$. If $\pi<\frac{1}{1+\gamma}$, then the $X$
fluctuations are of order $1/2$, and likewise if
$\eta<\frac{\gamma}{1+\gamma}$ then the $Y$ fluctuations are of
order $1/2$, whereas in the complementary cases, the fluctuations
are of order $1/3$. By comparing leading (law of large number) terms
in Propositions \ref{X_dist} and~\ref{Y_dist} we see that if either
of these two inequalities hold, then the fluctuations must be of
order $1/2$. In this case, either the leading term for $X$ or $Y$
clearly wins, in which case $L_2(N,M)$ has the leading order
behavior and Gaussian fluctuations of the winner random variable, or
the two random variables have the same leading terms. The second
case, or equal leading terms, occurs when
%
%
\begin{equation}\label{same_eqn}
\biggl(\pi^{-1}+\frac{\pi^{-1}\gamma^{2}}{\pi^{-1}-1} \biggr)N= \biggl(\eta
^{-1}+\frac{\eta^{-1}\gamma^{-2}}{\eta^{-1}-1} \biggr)M.
\end{equation}
In this case, the fluctuations will remain of order $1/2$, but will
behave as the fluctuations of two independent normal random
variables (what we call $G^2$).

If $\gamma=1$, there are two solutions to (\ref{same_eqn}).
One is $\eta=\pi$ and the other is $\eta=1-\pi$. Since we are only
considering the Gaussian region, the anti-diagonal solution is of no
interest, and we find that we have $G^2$ density for our
fluctuations if $\eta=\pi$ and $\eta<1/2$.

For $\gamma\neq1$, the solution set is a little harder. Recall
$M=\gamma^2 N$ and using this we can factor out $N$ from both sides
giving
%
%
\begin{equation}
\biggl(\pi^{-1}+\frac{\pi^{-1}\gamma^{2}}{\pi^{-1}-1} \biggr)= \biggl(\gamma
^2\eta^{-1}+\frac{\eta^{-1}}{\eta^{-1}-1} \biggr).
\end{equation}

Applying the change of variable $\pi\ra1-\pi$, we find that it
suffices to solve
%
%
\begin{equation}
\biggl(\gamma^2\pi^{-1}+\frac{\pi^{-1}}{\pi^{-1}-1} \biggr)= \biggl(\gamma
^2\eta^{-1}+\frac{\eta^{-1}}{\eta^{-1}-1} \biggr)
\end{equation}
and change the solution back to our original variables.

This again has the solution $\eta=\pi$. Solving for the other
solution and then changing variables back we get
%
%
\begin{equation}
\eta= \frac{\pi}{\pi(1-\gamma^{-2})+\gamma^{-2}}.
\end{equation}

By plugging in the critical point $(1/(1+\gamma),\gamma/(1+\gamma))$
it is easy to see that this $G^2$ curve is continuous between the
origin and the critical point, though only linear for $\gamma=1$.

We have, so far, only accounted for the regions where
$\pi<\frac{1}{1+\gamma}$ or $\eta<\frac{\gamma}{1+\gamma}$. There
are four other regions to consider which correspond to replacing the
or with an and, and the less than sign with either equality, or a
greater than sign. In each of these cases, the leading term is
independent of $\pi$ and $\eta$ and equals $(1+\gamma^{-1})^2M$. The
fluctuations of $X$ and $Y$ are both or order $1/3$, so those of
$L_2(N,M)$ are as well. To determine the scaling functions, finer
coupling arguments are necessary. For instance, when
$\pi>\frac{1}{1+\gamma}$ and $\eta>\frac{\gamma}{1+\gamma}$, the
last passage path for $X$ and for $Y$ can be compared to the
analogous random variables $\tilde{X}$ and $\tilde{Y}$, for last
passage paths for a coupled LPP model with two-sided boundary
conditions all of rate one (i.e., $\pi,\eta=1$) such that pointwise
$X\geq\tilde{X}$ and likewise $Y\leq\tilde{Y}$. The maximum of
$\tilde{X}$ and $\tilde{Y}$ equals $\tilde{L}(N,M)$, the last
passage time, and we show that $X$ and $\tilde{X}$, and likewise for
$Y$ and $\tilde{Y}$, under appropriate centering and scaling
converge to the same distribution, respectively. This implies, via
our Lemma \ref{max_lemma}, that the scaled and centered random
variables in fact converge in probability and hence that their
maximums converge in probability. This means that the maximum of $X$
and $Y$ behaves just like a regular last passage time which is known
to have $\Gue$ scaling function. A similar coupling shows that the
scaling when either $\pi=\frac{1}{1+\gamma}$ and
$\eta>\frac{\gamma}{1+\gamma}$, or $\pi>\frac{1}{1+\gamma}$ and
$\eta=\frac{\gamma}{1+\gamma}$ the scaling\vspace*{1pt} function behaves like
that of a single critical last passage time for a LPP model with
one-sided boundary conditions.

Determining the scaling function at the critical point is a harder
problem. One may identify it as the maximum of two $F_1$
distributions, coupled as $X$ and $Y$ are coupled. This
characterization, a priori, yields a tight family of random
variables. However, how to prove that it converges on more than just
a subsequence is not immediately clear, and more over it is not
clear to what it convergences. A posteriori, this characterization is
justified since the result of \cite{PLFS2006s} can readily be
translated into a proof that the scaling function at the critical
point is $F_{1,1}(x;0;0)$ (what they call $F_0$).

The results from \cite{BBP2005p} used in the proof of Theorem
\ref{LPP_fluctuations} yield, in fact, a much more general result
via essentially the same argument. We now define what we call the
\textit{LPP model with thick two-sided boundary conditions} in terms
of boundary row and column thickness integer parameters $J,I\geq0$;
two vectors of row and column weight rates $\pi=(\pi_1,\ldots,
\pi_J)$, $\eta=(\eta_1,\ldots,\eta_I)$; two vectors of row and
column convergence rates $X=(x_1,\ldots, x_J)$, $Y=(y_1,\ldots,
y_I)$. With these parameters, our model is defined in terms of the
following LPP weights (which implicitly depend on $M$ and $N$):
%
%
\begin{equation} \label{LPP_weights_thick}
w_{i,j} = \cases{
\mbox{exponential of rate }\pi_j + \dfrac{x_j}{M^{1/3}}, &\quad if
$i>I,j\leq J$;
\vspace*{2pt}\cr
\mbox{exponential of rate }\eta_i + \dfrac{y_i}{M^{1/3}}, &\quad if $i\leq
I,j> J$;
\vspace*{2pt}\cr
\mbox{exponential of rate }1, &\quad if $i> I,j> J$; \vspace*{2pt}\cr
\mbox{zero}, &\quad if $i\leq I,j\leq J$.}
\end{equation}

To see that this model is a broad generalization of our previously
considered two-sided boundary condition model, take $J,I=1$,
$\pi_1=\pi$, $\eta_1=\eta$ and $x_1$, \mbox{$y_1=0$}. Corresponding to this
model, we now provide a complete characterization of its asymptotic
fluctuations. A number of distributions not previously discussed are
introduced in this theorem. A full discussion of these distributions
can be found in Section \ref{notation}.

The coupling arguments given to prove Theorem \ref{LPP_fluctuations}
can be easily adopted to this new setting. Given the above
parameters define, again, two coupled random variables $X(N,M)$ and
$Y(N,M)$ as follows. $X(N,M)$ is the last passage time from $(0,0)$
to $(N,M)$ of the set of up/right paths which cross through at least
one vertex from the set $\{(i,j)\dvtx i=I, j\in\{0,\ldots,J-1\}\}$.
Likewise $Y(N,M)$ is the last passage time from $(0,0)$ to $(N,M)$
of the set of up/right paths which cross through at least one vertex
from the set $\{(i,j)\dvtx i\in\{0,\ldots,I-1\}, j=J\}$. Clearly, any
up/right path from $(0,0)$ to $(N,M)$ must go through one and only
one of these two regions. Furthermore, by virtue of the definition
of the last passage time, the maximizing path for $X(N,M)$ and
$Y(N,M)$ will necessarily go through the points $(0,J)$ and $(I,0)$
(resp.). Thus we may refine our definitions of $X(N,M)$ and
$Y(N,M)$ to require passing through these two points. Again, we see
that $L_2(N,M)=\max(X(N,M),Y(N,M))$ and just as before
\cite{BBP2005p} provides an immediate proof of the following.
\begin{proposition}\label{X_dist_thick}
For the vector $\pi$, fix the set $K_1\subset\{1,\ldots,J\}$ by
%
%
\begin{equation}
K_1=\biggl\{j\in\{1,\ldots,J\}\dvtx\pi_{j}=\frac{1}{1+\gamma}\biggr\}
\end{equation}
an define $X_{K_1}$ as the elements of $X$ which correspond to
indices in $K_1$. Further, define
$\tilde{\pi}=\min_{j\in\{1,\ldots,J\}}(\pi_j)$ and let $k_1$ be the
number of $\pi_j$ which attain the value~$\tilde{\pi}$.

Then with $X(N,M)$ defined as above, as $M,N\ra\infty$ while
$M/N-\gamma^2$ is in a compact subset of $(0,\infty)$, the following
holds for any real $x$ in a compact set:
\begin{enumerate}[(2)]
\item[(1)] when $\tilde{\pi}>\frac{1}{1+\gamma}$,
%
%
\begin{equation}
P \biggl(X(N,M)\leq(1+\gamma)^2N + \frac{(1+\gamma^{-1})^{4/3}}{\gamma
^{-1}}N^{1/3} x \biggr) \ra F_0(x);
\end{equation}
\item[(2)] when $\tilde{\pi}=\frac{1}{1+\gamma}$,
%
%
\begin{equation}
P \biggl(X(N,M)\leq(1+\gamma)^2N + \frac{(1+\gamma^{-1})^{4/3}}{\gamma
^{-1}}N^{1/3} x \biggr) \ra F_{|K_1|}(x;X_{K_1});
\end{equation}
\item[(3)] when $\tilde{\pi}<\frac{1}{1+\gamma}$,
%
%
\begin{eqnarray}
&&
P \biggl(X(N,M)\leq
\biggl(\tilde{\pi}^{-1}+\frac{\tilde{\pi}^{-1}\gamma^{2}}{\tilde{\pi
}^{-1}-1} \biggr)N
+
\biggl(\tilde{\pi}^{-2}-\frac{\tilde{\pi}^{-2}\gamma^{2}}{(\tilde{\pi
}^{-1}-1)^2} \biggr)^{1/2}
N^{1/2} x \biggr)\nonumber\\[-8pt]\\[-8pt]
&&\qquad\ra G_{k_1}(x).\nonumber
\end{eqnarray}
\end{enumerate}
\end{proposition}

A similar proposition exists for $Y(N,M)$. Using these two results,
the same types of coupling arguments then apply and give both the
orders and the scaling functions for $L_2(N,M)$. As before, these
coupling arguments break down when both boundary conditions are
critical. With single width boundary conditions, we appealed to
\cite{PLFS2006s}, however in this case no existing argument
provides a characterization of the behavior in this case. The
following partial theorem therefore contains a single conjectured
equation (\ref{partial_res}) whose study seems very difficult.
\begin{ptheorem}\label{LPP_fluctuations_thick}
\begin{enumerate}[(2)]
\item[(1)]
For vectors $\pi,\eta$ fix the sets $K_1\subset\{1,\ldots, J\}$ and
$K_2\subset\{1,\ldots, I\}$ by
%
%
\begin{eqnarray}
K_1&=&\biggl\{j\in\{1,\ldots,J\}\dvtx\pi_j=\frac{1}{1+\gamma}\biggr\},\\
K_2&=&\biggl\{i\in\{1,\ldots,I\}\dvtx\eta_i=\frac{\gamma}{1+\gamma}\biggr\}.
\end{eqnarray}
Then define $X_{K_1}$ and $Y_{K_2}$ as the elements of $X$ and $Y$
which correspond to indices in $K_1$ and $K_2$, respectively.

For $\gamma\in(0,\infty)$ and $M/N\ra\gamma^2$, then for vectors
$\pi,\eta$ such that $\pi_j\geq\frac{1}{1+\gamma}$ for all
$i\in\{1,\ldots,J\}$ and $\eta_i\geq\frac{\gamma}{1+\gamma}$ for all
$i\in\{1,\ldots,I\}$ then if:

\begin{enumerate}[(a)]
\item[(a)] $|K_1|=0$, $|K_2|=0$,
%
%
\begin{equation}
P\biggl(L_2(N,M)\leq(1+\gamma^{-1})^2M +\frac{(1+\gamma)
^{4/3}}{\gamma}M^{1/3}x\biggr)\ra F_0(x);
\end{equation}
\item[(b)] $|K_1|>0$, $|K_2|=0$,
%
%
\begin{equation}
P\biggl(L_2(N,M)\leq(1+\gamma^{-1})^2M
+\frac{(1+\gamma)^{4/3}}{\gamma}M^{1/3}x\biggr)\ra F_{|K_1|}(x;X_{K_1});
\end{equation}
\item[(c)] $|K_1|=0$, $|K_2|>0$,
%
%
\begin{equation}
P\biggl(L_2(N,M)\leq(1+\gamma^{-1})^2M
+\frac{(1+\gamma)^{4/3}}{\gamma}M^{1/3}x\biggr)\ra F_{|K_2|}(x;Y_{K_2});
\end{equation}
\item[(d)] $|K_1|>0$, $|K_2|>0$,
%
%
\begin{eqnarray}\label{partial_res}
&&P\biggl(L_2(N,M)\leq(1+\gamma^{-1})^2M
+\frac{(1+\gamma)^{4/3}}{\gamma}M^{1/3}x\biggr)\nonumber\\[-8pt]\\[-8pt]
&&\qquad\ra
F_{|K_1|,|K_2|}(x;X_{K_1},Y_{K_1}).\nonumber
\end{eqnarray}

\end{enumerate}
\item[(2)]
Define $\tilde{\pi} = \min_{j\in\{1,\ldots,J\}}(\pi_j)$ and
$\tilde{\eta} = \min_{i\in\{1,\ldots,I\}}(\eta_i)$, and let $k_1$
be the number of $\pi_j$ which attain the value $\tilde{\pi}$ and
likewise $k_2$ be the number of $\eta_i$ which attain the value
$\tilde{\eta}$.

For $\gamma\in(0,\infty)$ and $M/N\ra\gamma^2$, then for
$\pi,\eta$ such that:

\begin{enumerate}[(a)]
\item[(a)] $\tilde{\pi}<1/(1+\gamma)$ and $\tilde{\eta}>\frac{\tilde
{\pi
}}{\tilde{\pi}(1-\gamma^{-2})+\gamma^{-2}}$ [the $G$ ($\pi$
controlled) region],
%
%
\begin{eqnarray}
&&P \biggl(L_2(N,M)\leq
\biggl(\tilde{\pi}^{-1}+\frac{\tilde{\pi}^{-1}\gamma^{2}}{\tilde{\pi
}^{-1}-1} \biggr)N\nonumber\\
&&\qquad\hspace*{51.4pt}{}+
\biggl(\tilde{\pi}^{-2}-\frac{\tilde{\pi}^{-2}\gamma^{2}}{(\tilde{\pi
}^{-1}-1)^2} \biggr)^{1/2}
N^{1/2} x \biggr)\\
&&\qquad\ra G_{k_1}(x);\nonumber
\end{eqnarray}

\item[(b)] $\tilde{\eta}<\gamma/(1+\gamma)$ and $\tilde{\eta}<\frac
{\tilde
{\pi}}{\tilde{\pi}(1-\gamma^{-2})+\gamma^{-2}}$ [the $G$ ($\tilde
{\eta
}$ controlled) region],
%
%
\begin{eqnarray}
&&P \biggl(L_2(N,M)\leq
\biggl(\tilde{\eta}^{-1}+\frac{\tilde{\eta}^{-1}\gamma^{-2}}{\tilde
{\eta
}^{-1}-1} \biggr)M\nonumber\\
&&\qquad\hspace*{51pt}{}+
\biggl(\tilde{\eta}^{-2}-\frac{\tilde{\eta}^{-2}\gamma^{-2}}{(\tilde
{\eta
}^{-1}-1)^2} \biggr)^{1/2}
M^{1/2} x \biggr)\\
&&\qquad\ra G_{k_2}(x);\nonumber
\end{eqnarray}
\item[(c)] $\tilde{\pi}+\tilde{\eta}<1$ and $\tilde{\eta}=\frac
{\tilde{\pi
}}{\tilde{\pi}(1-\gamma^{-2})+\gamma^{-2}}$ (the $G^2$ line),
%
%
\begin{eqnarray}
&&P \biggl( L_2(N,M)\leq\biggl(\tilde{\pi}^{-1}+\frac{\tilde{\pi}^{-1}\gamma
^{2}}{\tilde{\pi}^{-1}-1} \biggr)N \nonumber\\[-8pt]\\[-8pt]
&&\qquad\hspace*{51.5pt}{} + \biggl(\frac{(1-\tilde{\pi} +\tilde
{\pi}\gamma^2)((1-\tilde{\pi})^2-\tilde{\pi}^2\gamma^2)}{\gamma
^2\tilde
{\pi}^2(1-\tilde{\pi})^2} \biggr)^{1/2} N^{1/2} x \biggr)\nonumber\\
&&\qquad\ra G_{k_1}\biggl(x\frac{\gamma}{\sqrt{1-\tilde{\pi}+\gamma^2
\tilde{\pi}}}\biggr)G_{k_2}\biggl(x \frac{\sqrt{1-\tilde{\pi}+\gamma^2
\tilde{\pi}}}{\gamma}\biggr).
\end{eqnarray}
\end{enumerate}
\end{enumerate}
\end{ptheorem}

Finally, let us note two applications of LPP with thick one-sided
boundary conditions which can be found in \cite{JB2006p}. The first
application deals with what Baik called \textit{traffic of slow start
from stop} in which particles start in the step initial condition of
TASEP and have a start-up profile---that is, every particle moves
slower for its first few jumps, and then returns to jumping at rate
one. The second application is dual to the first one and is called
\textit{traffic with a few slow cars} in which particles always move
at a slower rate. In both cases, Baik identifies the fluctuation
scaling limits by using the \cite{BBP2005p} type results which we
have made use of herein.

In the next section, we will give an important application for the
LPP with two-sided boundary conditions model to two-sided TASEP. It
is unclear whether the thick two-sided boundary conditions model has
any similar application to TASEP or related models.

\section{Mapping TASEP to last passage percolation with boundary
conditions}\label{Mapping_section}
In this section, we explain the connections between the last passage
time in LPP with two-sided boundary conditions and the fluctuations
of the height function for the
two-sided TASEP model. Making use of this mapping, we explain how the
results of Theorem \ref{LPP_fluctuations} imply the results of
\cite{PAFF1994c} stated in the \hyperref[sec1]{Introduction}.
Furthermore, we briefly
explain how this theorem translates into a proof of Conjecture 7.1
of \cite{PS2002c} (full proof is given in Section \ref{proofs}).

We start with a lemma which states that finite perturbations of our
LPP model, have no affect on the asymptotic behavior of the last
passage time.
\begin{lemma}\label{finite_weight_changes}
Fix some LPP model with weights $w_{i,j}$ (independent but not
necessarily identically distributed) such that
%
%
\begin{equation}
P \biggl(\frac{L(N,M)-a_N}{b_N}\leq x \biggr)\ra F(x)
\end{equation}
for $M/N\ra\gamma^2\in(0,\infty)$, and for $F$ a nondegenerate
probability distribution. Randomly, independent of the values of
$w_{i,j}$, change a set of these weights to a new set of weights
$w'_{i,j}$ and let $L'(N,M)$ denote the last passage time with
respect to the original weight with the newly updated weights. Call
$A$ the set of changed indices $(i,j)$ and $W_A = \sum_{(i,j)\in A}
w_{i,j}+w'_{i,j}$. Then if $E[W_A]<\infty$ and if $b_N\ra\infty$,
we also have
%
%
\begin{equation}
P \biggl(\frac{L'(N,M)-a_N}{b_N}\leq x \biggr)\ra F(x).
\end{equation}
\end{lemma}

Below is an outline of the proof. The full level of details is
suppressed since a similar style of proof is given for Lemma
\ref{stoch_dom_convergence} in full detail.
\begin{pf*}{Proof of Lemma \ref{finite_weight_changes}}
Since the total effect of the change of weights corresponding to $A$
has finite expectations, the Markov inequality shows that for any $\e$
we can find $l$ large enough so that $P(W_A\geq l)\leq\e$. If we
restrict ourselves to this region of our statespace, then since the
$b_N$ goes to infinity, the effect of the change of weights is
negligible in the limit. Since this is true on all but an $\e$ region
of the state space, we have that the distribution functions are within
$\e$ of each other in the limit, but taking $\e$ to zero gives equality.
\end{pf*}

Recall our definition of the two-sided TASEP model given by the
initial conditions of Bernoulli with parameter $\rho_-$ on the left
of zero and with parameter $\rho_+$ on the right. Corresponding to the
TASEP process started with
this random initial condition, we consider the height function
$h_t(j)$ defined in (\ref{height_function}). Theorem 2.1 of
\cite{PS2002c} relates the joint distributions for this height
function to those of the height function for a particular growth
model associated with a variant on the LPP with two-sided boundary
conditions. The weights for this variant LPP are defined with
respect to two independent geometric random variables $\zeta_+$ and
$\zeta_-$. Let $\zeta_+$ be geometric with parameter $1-\rho_+$
[i.e., $P(\zeta_+=n)=\rho_+(1-\rho_+)^n$] and $\zeta_-$ be geometric
with parameter $\rho_-$ [i.e., $P(\zeta_-=n)=(1-\rho_-)\rho_+^n$].
The weights are then defined as independent random variables with:
%
%
\begin{equation}\label{fullbcs}
w_{i,j}=\cases{
\mbox{exponential of rate 1}, &\quad if $i,j\geq1$;\cr
\mbox{zero}, &\quad if $i=j=0$;\cr
\mbox{zero}, &\quad if $0\leq i\leq\zeta_{+}\mbox{ and } j=0$;\cr
\mbox{exponential of rate }1-\rho_+, &\quad if $i>\zeta_+\mbox{ and }
j=0$;\cr
\mbox{zero}, &\quad if $0\leq j\leq\zeta_-\mbox{ and } i=0$;\cr
\mbox{exponential of rate }\rho_-, &\quad if $j>\zeta_-\mbox{ and } i=0$.}
\end{equation}

With respect to these random weights, define the last passage time
$\tilde{L}(N,M)$. This family of random variables is nondecreasing in
both $N$ and $M$. Therefore, one can associate to this a growth
process on the upper-corner and likewise a height process over the number line.
Let $A_t = \{(N,M)|N,M\geq1,\tilde{L}(N,M)\}$ be the growth process\vspace*{1pt}
and let $\tilde{h}_t$ be defined so as to satisfy $A_t=\{(N,M)|2\leq
N+M\leq\tilde{h}_t(N-M)\}$. To describe this in words, imagine
rotating counter-clockwise, the upper corner in which LPP occurs by
$\pi
/4$. To each lattice point (labeled by $j\in\Z$) now on the horizontal
associate a height $\tilde{h}_t(j)$ equal to two times the number of
$\tilde{L}(N,M)$ vertically above $j$ which are less than or equal to
$t$. Then we have the following.
\begin{theorem}[(Theorem 2.1 of \cite{PS2002c})] \label{Joint_dist}
In the sense of joint distributions, we have
%
%
\begin{equation}
h_{t}(j)=\tilde{h}_t(j) \qquad\mbox{for } |j|\leq h_{t}(j).
\end{equation}
\end{theorem}

This theorem essentially says that for the height profile which lies
above the boundary of the rotated upper corner, the two profiles
have the same joint distribution. It is worthwhile to recall that
there is a similar map between the TASEP height function for TASEP
with step initial conditions and the height function for standard
(no boundary condition) LPP \cite{KJ2000s}. The proof of Theorem
\ref{Joint_dist} can be found in \cite{PS2002c} and essentially
amounts to a study of the dynamics of the right most particle to the
left of the origin, as well as the dynamics of the left most hole to
the right of the origin. Tagging this particle and this hole, we
observe that their initial location is geometric and using Burke's
theorem we find that their waiting times between successive moves is
exponential with rate relating to the densities $\rho_-$ and~$\rho_+$.
Between the tagged particle on the left and the tagged
hole on the right, particles move according to normal TASEP rules,
and hence the two height functions evolve with the same dynamics.
The dynamics of the boundary of the part of the TASEP height
function lying in the rotated upper corner is matched by the effect
of the LPP boundary conditions, and the theorem follows.

From this theorem, we see that the following equality:
%
%
\begin{equation}\label{height_to_LPP}
P_{\rho_-,\rho_+}\bigl(h_{t}(N-M)\geq N+M\bigr) = P\bigl(\tilde{L}(N,M)\leq t\bigr).
\end{equation}

This equality is not of much use to use, however, because in order to
use the results of Theorems \ref{LPP_fluctuations}, we must deal with
a slightly different LPP model than corresponding to $\tilde{L}$.
However, this model and our LPP with two-sided boundary conditions
only differ in expectation by a finite number of weights (the
geometric number of zeros from $\zeta_-$ and $\zeta_+$). Therefore,
while it is true that there is not exact equality then with
$P(L_2(N,M)\leq t)$, from Lemma \ref{finite_weight_changes}, we see
that for any sort of central limit fluctuation statement with a
nontrivial limiting distribution, we have equality in the limit. We
will abuse notation in the remainder of this section and the next
during the proof of Conjecture 7.1, and write equality between the
TASEP height function probability and the probability for the last
passage time in the LPP with two-sided boundary conditions model. To
sum up, we have the following.
\begin{rem}
While the boundary conditions (\ref{LPP_weights}) differ from those
(\ref{fullbcs}) used by Pr\"{a}hofer and Spohn, they are much simpler
and also describe the TASEP with two-sided initial conditions, as
described in \cite{BCS2006c}.
\end{rem}

It is important to note that Lemma \ref{finite_weight_changes} only
applies if both $\rho_-<1$ and $\rho_+>0$. If either of these
inequalities is violated, then the geometric number of zeros on the
boundary will in fact, almost surely be infinite. However, in any of
these cases, the classification of one-sided LPP then readily
applies.

The boundary conditions for $\tilde{L}$ corresponded to having
exponentials of rate $\rho_-$ on the left boundary and $1-\rho_+$ on
the right. Therefore, in terms of $\pi$ and $\eta$, we have
$\pi=1-\rho_+$ and $\eta=\rho_-$. The critical point for $\pi,\eta$
is $\frac{1}{1+\gamma}$ and $\frac{\gamma}{1+\gamma}$, therefore we
see that the critical point for $\rho_-,\rho_+$ is
%
%
\begin{equation}
\rho_-=\rho_+=\frac{\gamma}{1+\gamma}.
\end{equation}
This corresponds to an equilibrium measure on TASEP with density
$\frac{\gamma}{1+\gamma}$.

In the next section, we will show how Theorem \ref{LPP_fluctuations}
implies an almost complete (all but a few regions of the claimed
Gaussian region are fully proved) proof of~\cite{PS2002c},
Conjecture 7.1 (our Theorem \ref{thm_conj}). From this result, we may
easily deduce the results of \cite{PAFF1994c} stated in the
\hyperref[sec1]{Introduction}.

For simplicity assume $r\in[0,1]$, as the case $r\in[-1,0]$ follows
similarly. We wish to prove that
%
%
\begin{equation} \label{FF}
P \biggl(\frac{J_{rt,t}-(\rho(1-\rho)-r\rho)t}{t^{1/2}\sqrt{\rho
(1-\rho
)|(1-2\rho)-r|}} \geq x \biggr)\ra G_1(x),
\end{equation}
where as defined before $G_1(x)$ is the standard Gaussian
distribution function.

In the case of $\rho_-=\rho_+=\rho$ and $r\geq0$, we can conclude
from (\ref{conj_eq1}) that
%
%
\begin{equation}
P \bigl(t\bar{h}(r)-h_t([ry])\leq\bigl(4\rho(1-\rho)(r-1+2\rho)t\bigr)^{1/2}x
\bigr)\ra G_1(x),
\end{equation}
where $\bar{h}(r) = (1-2\rho)r + 2\rho(1-\rho)$. Substituting the
relationship in (\ref{height_to_current}) and rearranging
terms, we arrive at the exact result of \cite{PAFF1994c} desired.

We now briefly explain the approach to proving Conjecture 7.1 from
our Theorem \ref{LPP_fluctuations}. The conjecture deals with height
functions. We have provided above the relationship between height
function distributions and LPP distributions.
From~(\ref{height_to_LPP}), we see that if one is to consider
$P(h_t(j)\geq x)$ in terms of LPP, you must solve for $N =
\frac{x+j}{2}$ and $M=\frac{x-j}{2}$. The variables $j$ and $x$ both
are functions of time $t$ and a speed $y$. If $M/N =
\frac{x-j}{x+j}$ has a\vspace*{1pt} limit, we call that $\gamma^2$. This allows
us to asymptotically write $M$ (or $N$) just as a function of time
(thus the $y$ dependence goes into $\gamma$). In the cases we
consider, we can invert the expression for $M$ in terms of $t$ and
get an expression for $t$ in terms of $M$, thus putting us in the
form of the limit theorems we proved in Theorem
\ref{LPP_fluctuations}.

\section{Proof of fluctuation theorems}\label{proofs}
In this section, we provide a proof of Theorem \ref{LPP_fluctuations}
(which easily generalizes to prove Partial Theorem
\ref{LPP_fluctuations_thick}) and a proof of Theorem \ref{thm_conj}.

\subsection[Proof of Theorem 1.3
(fluctuations of LPP with two-sided boundary conditions]{Proof of
Theorem \protect\ref{LPP_fluctuations}
(fluctuations of LPP with two-sided boundary conditions)}

The following two technical lemmas provide the basis for the
coupling arguments necessary in our proof of Theorem
\ref{LPP_fluctuations}.
\begin{lemma}\label{stoch_dom_convergence}
If $X_n\geq\tilde{X}_n$ and $X_n\Rightarrow D$ as well as $\tilde
{X}_n\Rightarrow D$, then $X_n-\tilde{X}_n$ converges to zero in
probability. Conversely if $X_n\geq\tilde{X}_n$ and $\tilde
{X}_n\Rightarrow D$ and $X_n-\tilde{X}_n$ converges to zero in
probability, then $X_n\Rightarrow D$ as well.
\end{lemma}
\begin{lemma}\label{max_lemma}
Assume $X_n\geq\tilde{X}_n$ and $X_n\Rightarrow D_1$ as well as
$\tilde{X}_n\Rightarrow D_1$; and similarly $Y_n\geq\tilde{Y}_n$ and
$Y_n\Rightarrow D_2$ as well as $\tilde{Y}_n\Rightarrow D_2$. Let
$Z_n=\max(X_n,Y_n)$ and $\tilde{Z}_n=\max(\tilde{X}_n,\tilde{Y}_n)$.
Then if $\tilde{Z}_n\Rightarrow D_3$, we also have $Z_n\Rightarrow D_3$.
\end{lemma}
\begin{pf*}{Proof of Lemma \ref{stoch_dom_convergence}}
While it is likely that this lemma is known in the literature, we do
not know where and hence produce a proof.
We prove the first assertion. Fix $\e>0$ and, from the point of
contradiction, assume that $P(X_n-\tilde{X}_n>\e)>\delta>0$ for an
infinite subsequence of $n$'s. By restricting to this subsequence
and noting that all of the hypothesis of the lemma hold under this
restriction, we may equivalently assume that
$P(X_n-\tilde{X}_n>\e)>\delta>0$ for all $n$ large. Since $X_n$
and $\tilde{X}_n$ converge weakly, each sequence of random variables
is tight. This implies that there exists an $M(\e)$ and $N(\e)$ such
that for all $n>N$, $P(|X_n|>M)<\delta/2$ and likewise
$P(|\tilde{X}_n|>M)<\delta/2$. Thus $P(|\tilde{X}_n|>M \cap
\{X_n-\tilde{X}_n>\e\})<\delta/2$, therefore
%
%
\begin{equation}P(|\tilde{X}_n|<M \cap
\{X_n-\tilde{X}_n>\e\})>\delta/2.
\end{equation}
Call this event $A=|\tilde{X}_n|<M \cap\{X_n-\tilde{X}_n>\e\}$,
then conditioned on $A$, $X_n>\tilde{X}_n+\e$. For large enough $n$,
%
%
\begin{equation}P(X_n\leq t|A)\leq P(\tilde{X}_n\leq t|A)-P(\tilde
{X}_n\in[t-\e,t]|A).
\end{equation}

We now partition the interval $[-M,M]$ into $\e$ size blocks and define
deterministic numbers $a_j(n)$ for $j\in\{1,\ldots, \lceil\frac
{2M}{\e
} \rceil\}$ by
%
%
\begin{equation}
a_j(n)=P\bigl(\tilde{X}_n\in[-M+\e(j-1),-M+\e j]|A\bigr).
\end{equation}
Observe that $\sum_{j} a_j(n)=1$ since having conditioned on $A$, we
know $\tilde{X}_n\in[-M,M]$. Therefore, for each $n$, there exists at
least one $j=j(n)$ for which $a_j(n)\geq
\frac{1}{2M/\e+1}=\frac{\e}{2M+\e}$ [if there is more than one $j$ for
which $a_j(n)$ is as desired, pick the smallest value of $j$]. Since
$j$ is restricted to a finite set of values, there must be some
infinite subsequence of $n$'s which have the same value of $j(n)$.
Restricting to that subsequence so every $j(n)$ equals a fixed $j$, if
we set $t=-M+\e j$ we have
%
%
\begin{equation}
P(\tilde{X}_n\in[t-\e,t]|A)\geq\frac{\e}{2M+\e}.
\end{equation}

Therefore,
%
%
\begin{equation}
P(X_n\leq t|A)\leq P(\tilde{X}_n\leq t|A)-\frac{\e
}{2M+\e}.
\end{equation}

Multiplying both sides by $P(A)$ and rewriting without
conditioning gives
%
%
\begin{equation}
P(X_n\leq t\cap A)\leq P(\tilde{X}_n\leq t \cap A) -
\frac{P(A)\e}{2M+\e}.
\end{equation}

That $X_n\geq\tilde{X}_n$ also implies that
%
%
\begin{equation}
P(X_n\leq t\cap A^c)\leq P(\tilde{X}_n\leq t \cap A^c).
\end{equation}

Adding these two inequalities and using the fact that
$P(A)>\delta/2$ gives, for all $n$ large enough
%
%
\begin{equation}
P(X_n\leq t)\leq P(\tilde{X}_n\leq t) -\frac{\delta\e
}{2(2M+\e)}.
\end{equation}

This inequality implies, however, that $X_n$ and $\tilde{X}_n$
cannot converge in distribution to the same object. This is a
contradiction to our hypothesis, so our assumption must be false.
That is, $P(X_n-\tilde{X}_n>\e)$ must go to zero as $n$ goes to
infinity.

The second assertion is easier. For all $\e$, we can find $N$ such
that for $n>N$, $P(X_n-\tilde{X}_n>\e)\leq\e$. Set inclusion and
partitioning implies that
%
%
\begin{eqnarray}
P(\tilde{X}_n\leq t-\e) &=& P(\tilde{X}_n\leq t-\e\cap X_n-\tilde
{X}_n<\e)\nonumber\\
&&{} + P(\tilde{X}_n\leq t-\e\cap X_n-\tilde{X}_n\geq\e
)\nonumber\\[-8pt]\\[-8pt]
&\leq& P(X_n\leq t)+P(X_n-\tilde{X}_n\geq\e)\nonumber\\
&\leq& P(X_n\leq t)+\e.\nonumber
\end{eqnarray}

Since $P(X_n\leq t)\leq P(\tilde{X}_n\leq t)$, we find that
%
%
\begin{equation}
P(\tilde{X}_n\leq t-\e)-\e\leq P(X_n\leq t) \leq P(\tilde{X}_n\leq t).
\end{equation}

If $t$ is any continuity point for $D$, then we can take $\e$ to
zero and we find that
%
%
\begin{equation}
D(t)\geq P(X_n\leq t) \geq D(t),
\end{equation}
and hence $X_n$ weakly converges to $D$.
\end{pf*}
\begin{pf*}{Proof of Lemma \ref{max_lemma}}
Applying Lemma \ref{stoch_dom_convergence} to both $X_N$ and $\tilde
{X}_N$, as well as $Y_N$ and $\tilde{Y}_N$ we find that for any $\e$,
large enough $N$, $P(A_X)<\e$ and likewise $P(A_Y)<\e$ where $A_X=\{
X_N-\tilde{X}_N>\e\}$ and $A_Y=\{Y_N-\tilde{Y}_N>\e\}$.
From this, it follows that
%
%
\begin{eqnarray}
&&P(\tilde{Z}_n\leq t-\e) \nonumber\\[-8pt]\\[-8pt]
&&\qquad= P(\tilde{Z}_n\leq t-\e\cap A_X^c\cap A_Y^c)+
P(\tilde{Z}_n\leq t-\e\cap A_X^c\cap A_Y)\nonumber\\
&&\qquad\quad{} + P(\tilde{Z}_n\leq t-\e\cap A_X\cap A_Y^c) + P(\tilde{Z}_n\leq t-\e
\cap A_X\cap A_Y).
\end{eqnarray}
The first probability is less than or equal to $P(Z_n\leq t)$ while the
last three are each trivially bounded by $\e$. Therefore, noting that
$P(Z_n\leq t)\leq P(\tilde{Z}_n\leq t)$ we find
%
%
\begin{equation}
P(\tilde{Z}_n\leq t-\e)-3\e\leq P(Z_n\leq t) \leq
P(\tilde{Z}_n\leq t).
\end{equation}

Taking $t$ to be a continuity point for the
limiting distribution for $\tilde{Z}_n$ (for $D_3$) and taking $\e$ to
zero we get that $\lim_{n\ra\infty}P(Z_n\leq t)=F_{D_3}(t)$, and hence
$Z_n$ converges in
distribution to $D_3$.
\end{pf*}
\begin{pf*}{Proof of Theorem \ref{LPP_fluctuations}, ($F_0$)}
This result follows from a coupling argument between the $X(N,M)$,
$Y(N,M)$ variables as well as a second set of last passage times
$\tilde
{X}(N,M)$, $\tilde{Y}(N,M)$. $X(N,M)$ and $Y(N,M)$ are coupled as
previously described (they are the last passage times if forced to go
right (or up) on the first move). Now to define the tilde versions of
$X(N,M)$ and $Y(N,M)$, for a given realization of weights, divide the
boundary weights by their means. This creates a new set of weights
coupled and pointwise dominated by the original set of weights. For
this new set of weights, define $\tilde{X}(N,M)$ and $\tilde{Y}(N,M)$
as the last passage time if forced right (or up) initially. From this
pointwise domination of the new weights by the original weights, we see
that $X(N,M)\geq\tilde{X}(N,M)$ and $Y(N,M)\geq\tilde{Y}(N,M)$
pointwise. The advantage of the tilde variables is that $\tilde
{Z}(N,M)=\max(\tilde{X}(N,M),\tilde{Y}(N,M))$ is the standard (without
boundary conditions) last passage time. We also know that,
asymptotically $X(N,M)$ and $\tilde{X}(N,M)$, as well as $Y(N,M)$ and
$\tilde{Y}(N,M)$ have the same distribution. We wish to use this
information to conclude that $Z(N,M)$ and $\tilde{Z}(N,M)$ have the
same distribution as well.

Let us redefine our variables by properly shifting and scaling them,
so that they have a limiting distribution. Our new $X(N,M)$ is
%
%
\begin{equation}
\frac{X(N,M)-(1+\gamma^2)N}{\gamma(1+\gamma^{-1})^{4/3}N^{1/3}},
\end{equation}
and similarly we define $\tilde{X}(N,M), Y(N,M),\tilde{Y}(N,M)$ and
$\tilde{Z}(N,M)$ and $Z(N,\break M)$ is redefined in terms of the newly
defined variables. Now we have the following setup: $X(N,M)\geq\tilde
{X}(N,M)$ and $X(N,M)\Rightarrow F_{0}$ as well as $\tilde
{X}(N,M)\Rightarrow F_{0}$; similarly $Y(N,M)\geq\tilde{Y}(N,M)$ and
$Y(N,M)\Rightarrow F_{0}$ as well as $\tilde{Y}(N,M)\Rightarrow F_{0}$.
Applying Lemma \ref{max_lemma}, we get that $Z(N,M)$ converges in
distribution to $F_{0}$.
\end{pf*}
\begin{pf*}{Proof of Theorem \ref{LPP_fluctuations}, ($F_1$)}
There are two cases which yield to the same argument. Thus, we prove
the case of $\pi>1/(1+\gamma)$ and $\eta=\gamma/(1+\gamma)$ only. As
in the last proof, let $X(N,M)$, $\tilde{X}(N,M)$, $Y(N,M)$ denote
the suitably centered and scaled random variables. For the sake of
applying Lemma \ref{max_lemma} we define $\tilde{Y}(N,M)=Y(N,M)$.
Again we know that $X(N,M)\geq\tilde{X}(N,M)$ and
$X(N,M)\Rightarrow F_{0}$ as well as $\tilde{X}(N,M)\Rightarrow
F_{0}$ and clearly the same holds for the $Y(N,M)$ and
$\tilde{Y}(N,M)$. So by Lemma \ref{max_lemma} since we know that
$\max(\tilde{X}(N,M),\tilde{Y}(N,M))$ converges weakly to $F_{1}$,
it follows that $\max(X(N,M),Y(N,M))\Rightarrow F_{1}$.
\end{pf*}
\begin{pf*}{Proof of Theorem \ref{LPP_fluctuations}, ($F_{1,1}$)}
This follows immediately from the main result of \cite{PLFS2006s}.
\end{pf*}
\begin{pf*}{Proof of Theorem \ref{LPP_fluctuations}, ($G$)}
We prove the first case, when $\pi<1/(1+\gamma)$ and
$\eta>\frac{\pi}{\pi(1-\gamma^{-2})+\gamma^{-2}}$, as the other case
has the same proof. In this region of $\pi,\eta$, the leading order
term on the expression for $X(N,M)$ is larger than that of the
expression for $Y(N,M)$. If we renormalize both $X(N,M)$ and
$Y(N,M)$ by the leading order term for $X(N,M)$ and divide by its
fluctuation term, we find that $X(N,M)$ converges to a standard
normal. On the other hand, since the leading order term for $X(N,M)$
exceeds that for $Y(N,M)$, the renormalized $Y(N,M)$ converges to
negative infinity. This implies that $\max(X(N,M),Y(N,M))$ converges
in distribution to a standard normal, just like $X(N,M)$.
\end{pf*}
\begin{pf*}{Proof of Theorem \ref{LPP_fluctuations}, ($G^2$)}
We couple $X(N,M)$ with a random variable $\tilde{X}(N,M)$ where
$\tilde{X}(N,M)$ is the last passage time when forced to stay along
the bottom edge for at least a specific, deterministic fraction of
the path [we likewise define $\tilde{Y}(N,M)$]. Specifically we
define $\tilde{X}(N,M)$ to be the last passage time when the path is
pinned to the bottom edge until the point
%
%
\begin{equation}
\biggl(1-\frac{\gamma^2}{(\pi^{-1}-1)^2} \biggr) N,
\end{equation}
after which point is forced into the bulk and allowed to follow a
last passage path therein. The weights accrued along the bottom edge
respect a simple central limit (as they are the sums of a
deterministic number of i.i.d. random variables) and the weights
accrued after the path is forced into the bulk follows the
fluctuations theorem for standard exponential last passage times. As
these two random variables are independent by construction, their
means add since their fluctuations are of different order ($N^{1/2}$
for the bottom and $N^{1/3}$ for the bulk) the bottom fluctuations
win out. Following this idea, we find that the mean of
$\tilde{X}(N,M)$ is
$(\pi^{-1}+\frac{\pi^{-1}\gamma^2}{\pi^{-1}-1})N$ and the
fluctuations are normal with variance
%
%
\begin{equation}
\pi^{-2} \biggl(1-\frac{\gamma^2}{(\pi^{-1}-1)^2} \biggr) N.
\end{equation}

To see this, observe that if we center $\tilde{X}$ by the mean and
divide by the square root of the above variance, we are left with a
random variable of the form
$\tilde{X}_1(N,M)+N^{-1/6}\tilde{X}_2(N,M)$, where $\tilde{X}_1(N,M)$
converges to a normal, and $\tilde{X}_2(N,M)$ converges to a $\Gue$.
Since the second term has a prefactor which goes to zero, we see
that $P(\tilde{X}_1(N,M)+N^{-1/6}\tilde{X}_2(N,M)\leq l)$ can be
partitioned into a region of size $\e$ where $|\tilde{X}_2(N,M)|\geq
R$ and a region of size $1-\e$ where $|\tilde{X}_2(N,M)|< R$. On
the second region, we can replace\vspace*{1pt} $\tilde{X}_2(N,M)$ by $R$ and find
asymptotically that the probability differs from
$P(\tilde{X}(N,M)\leq l)$ by only $\e$. Taking $\e$ to zero gives
the desired convergence in distribution.

Therefore if we center $X(N,M)$ and $\tilde{X}(N,M)$ by the same
amount and renormalize by the same amount we get two random
variables which converge to the same distribution, despite the first
one being almost always larger than the second one. This is one of
the pieces we will need to apply Lemma \ref{max_lemma}.

We can likewise define $\tilde{Y}(N,M)$ as the last passage time
when the path is pinned to the left edge until the point
%
%
\begin{equation}
\biggl(1-\frac{\gamma^{-2}}{(\eta^{-1}-1)^2} \biggr) M,
\end{equation}
after which point is forced into the bulk and allowed to follow a last
passage path therein. As in\vspace*{1pt} the prior we see that centered and
renormalizing $Y(N,M)$ and $\tilde{Y}(N,M)$ by the same amounts, gives
to random variables which converge to the same distribution, despite
the first one being almost always larger than the second one.

From the relationship between $\pi$ and $\eta$ we know that the
leading order terms for both $\tilde{X}(N,M)$ and $\tilde{Y}(N,M)$
coincide. Therefore, we can write
%
%
\begin{eqnarray}
\tilde{X}(N,M) &=& AN+B{\pi}N^{1/2} \tilde{X}_1(N,M) +C_{\pi}N^{1/3}
\tilde{X}_2(N,M),\\
\tilde{Y}(N,M) &=& AN+B{\eta}N^{1/2} \tilde{Y}_1(N,M) +C_{\eta}N^{1/3}
\tilde{Y}_2(N,M),
\end{eqnarray}
where
%
%
\begin{eqnarray}
A &=& \biggl(\pi^{-1} + \frac{\pi^{-1}\gamma^2}{\pi^{-1}-1}\biggr),\\
B_{\pi} &=& \biggl(\pi^{-2}-\frac{\pi^{-2}\gamma^2}{(\pi^{-1}-1)^2}
\biggr)^{1/2},\\
B_{\eta} &=& \biggl(\eta^{-2}-\frac{\eta^{-2}\gamma^{-2}}{(\eta
^{-1}-1)^2} \biggr)^{1/2},
\end{eqnarray}
and $C_{\pi}$ and $C_{\eta}$ are constants (which will play no role
here). If we consider now $\tilde{Z}(N,M) =
\max(\tilde{X}(N,M),\tilde{Y}(N,M))$, we find that
\[
P \biggl(\frac{\tilde{Z}(N,M)-AN}{N^{1/2}\sqrt{B_{\pi}B_{\eta}}}\leq
x \biggr)=P(E_1 \mbox{ and } E_2),
\]
where $E_1$ and $E_2$ are, respectively, the events
%
%
\begin{eqnarray}
\sqrt{\frac{B_{\pi}}{B_{\eta}}}\tilde{X}_1(N,M)+C_{\pi}'
N^{-1/6}\tilde{X}_2(N,M)&\leq& x,
\\
\sqrt{\frac{B_{\eta}}{B_{\pi}}}\tilde{Y}_1(N,M)+C_{\eta
}' N^{-1/6}\tilde{Y}_2(N,M)&\leq& x.
\end{eqnarray}

As before, because of the $N^{-1/6}$ prefactor to the
$\tilde{X}_2(N,M)$ and $\tilde{Y}_2(N,M)$ terms, we can condition on
these terms being bounded by some large number $R$, and only cost
ourselves $\e$ of the sample space. Once we have conditioned on
these random variables being bounded by $R$, we can conclude that
their joint probability is bounded between the product
%
%
\begin{eqnarray}
&&
P \Biggl(\tilde{X}_1(N,M)\leq
\sqrt{\frac{B_{\eta}}{B_{\pi}}}x-C'_{\pi}N^{-1/6}R \Biggr)\nonumber\\[-8pt]\\[-8pt]
&&\qquad{}\times P \Biggl(\tilde
{Y}_1(N,M)
\leq
\sqrt{\frac{B_{\pi}}{B_{\eta}}}x-C'_{\eta}N^{-1/6}R \Biggr)\nonumber
\end{eqnarray}
and
%
%
\begin{eqnarray}
&&
P \Biggl(\tilde{X}_1(N,M)\leq\sqrt{\frac{B_{\eta}}{B_{\pi
}}}x+C'_{\pi}N^{-1/6}R \Biggr)\nonumber\\[-8pt]\\[-8pt]
&&\qquad{}\times P \Biggl(\tilde{Y}_1(N,M)\leq\sqrt{\frac
{B_{\pi}}{B_{\eta}}}x+C'_{\eta}N^{-1/6}R \Biggr).\nonumber
\end{eqnarray}

Taking $N$ to infinity gives $P(\tilde{X}_1(N,M)\leq x)P(\tilde
{Y}_1(N,M)\leq x)$ and taking $\e$ to zero, and using the central limit
theorem to show that $\tilde{X}_1(N,M)$ is standard normal, we find that
%
%
\begin{equation}\label{max_eqn1}
P \biggl(\frac{\tilde{Z}(N,M)-AN}{N^{1/2}\sqrt{B_{\pi}B_{\eta}}}\leq
x \biggr)=G_1 \Biggl(x\sqrt{\frac{B_{\eta}}{B_{\pi}}} \Biggr)G_1 \Biggl(x\sqrt
{\frac{B_{\pi}}{B_{\eta}}} \Biggr).
\end{equation}

Therefore, using the observations made at the beginning of the
proof, we can apply Lemma \ref{max_lemma} to conclude that the
probability distribution of $Z(N,M)$ centered and normalized as
above, converges to the same $G^2$ distribution [i.e., equation
(\ref{max_eqn1}) holds with $Z(N,M)$ in place of $\tilde{Z}(N,M)$].
Working out the coefficients $B_{\pi}$ and $B_{\eta}$ using the
relationship between $\pi$ and $\eta$, we have our desired result.
\end{pf*}

\subsection{\texorpdfstring{Proof of Theorem \protect\ref{thm_conj} (Conjecture 7.1 from
\protect\cite{PS2002c})}{Proof of Theorem 1.1 (Conjecture 7.1 from [22])}} \label{proof_conj_71}

The following elementary lemma will find repeated use in
what follows.
\begin{lemma} \label{inverting}
If $M=at+bt^{1/2}$ then for large $t$,
%
%
\begin{equation}
t= a^{-1}M -a^{-3/2}bM^{1/2} + o(M^{1/2}).
\end{equation}
Likewise if $M=at+bt^{1/3}$ then for large $t$,
%
%
\begin{equation}
t= a^{-1}M -a^{-4/3}bM^{1/3} + o(M^{1/3}).
\end{equation}
\end{lemma}
\begin{pf*}{Proof of Theorem \ref{thm_conj}}
We provide proofs of only the $G^2$ and $F_0$ cases of this theorem
as the $G$ case is analogous to $G^2$ and the $F_1$ case is
analogous to $F_0$. The $F_{1,1}$ case of the theorem already is
proved in \cite{PLFS2006s}. The proofs are based on the fact that if
the height at a given time value exceeds a point $(N,M)$, then the last
passage time of that point is less than the above time value. Then
Theorem~\ref{LPP_fluctuations} applies and gives asymptotic height
distribution results. As noted before, the mapping between the TASEP
height function and the last passage time for our model of LPP with
two-sided boundary conditions is not exact (as the exact LPP model
has geometric numbers of boundary zeros) however Lemma
\ref{finite_weight_changes} ensures that asymptotically all results
in our LPP model correspond to results for the two-sided TASEP
height function.


\textit{$G^2$ case}: recall that in LPP with two-sided boundary
conditions $\pi= 1-\rho_+$ and $\eta= \rho_-$. We presently assume
that $\eta+\pi<1$ ($\rho_-<\rho_+$) and $y=y_c$. Recalling that
$\bar{h}(y)=(1-2\eta)y+2\eta(1-\eta)$, we wish to determine the
asymptotic (large $t$) value of
%
%
\begin{equation}
P\bigl(h_t(yt)\geq\bigl((1-2\eta)y+2\eta(1-\eta)\bigr)t - \sqrt{(1-\pi
-\eta)}t^{1/2}x\bigr).
\end{equation}
We will reduce this probability to a probability in the related LPP
with two-sided boundary conditions, and then use Theorem \ref
{LPP_fluctuations} to conclude that this probability is the correct
product of Gaussian probability functions.

%
%
\begin{figure}

\includegraphics{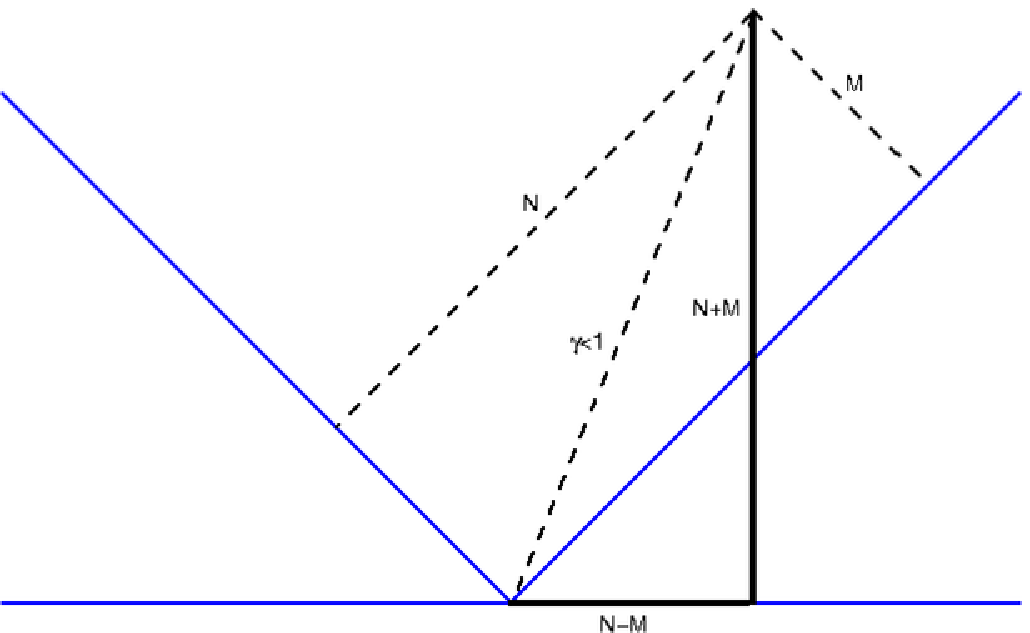}

\caption{The event that the last passage time $L(N,M)\leq t$
corresponds to the event that the height function $h_t(N-M)\geq
N+M$.}\label{tilted}
\end{figure}

The first step in translating to a LPP problem is to relate the
speed to $\gamma$. We may solve for the asymptotic value of $\gamma$
as a function of the speed $y$. As shown in Figure~\ref{tilted}, the
height function event corresponds to the LPP event where $N-M = yt$
and $N+M = ((1-2\eta)y+2\eta(1-\eta))t -
\sqrt{(1-\pi-\eta)}t^{1/2}x)$. From that, we find that
%
%
\begin{eqnarray}
M &=& \eta(1-y-\eta) t -\tfrac{1}{2}\sqrt{(1-\pi-\eta
)}t^{1/2} x,
\\
N &=& \bigl(y-\eta y +\eta(1-\eta)\bigr) t -\tfrac{1}{2}\sqrt{(1-\pi
-\eta)}t^{1/2} x.
\end{eqnarray}

Therefore,
%
%
\begin{equation}
\gamma^2 = \lim_{t\ra\infty} M/N = \frac{\eta(1-y-\eta
)}{y-\eta y +\eta(1-\eta)}.
\end{equation}

From this equation, we can solve for $y$ as a function of $\gamma$:
%
%
\begin{equation}
y = \frac{(\gamma^2-1)\eta(1-\eta)}{\gamma^2(\eta
-1)-\eta}.
\end{equation}

Since we have assume that $y=y_c$, we may use these two expressions
for $y$ to relate $\eta,\pi$ and $\gamma$ to find that
\begin{equation}
\eta= \frac{\pi}{\pi(1-\gamma^{-2})+\gamma^{-2}}.
\end{equation}
This is exactly the curve along which the $G^2$ part of Theorem \ref
{LPP_fluctuations} applies.

Finally, we may use Lemma \ref{inverting} to invert our expression
for $N$ in terms of $t$. Asymptotically
%
%
\begin{equation}\quad
t = \frac{1+\pi(\gamma^2-1)}{\pi(1\pi)}N + \frac{(1-\pi
-\eta)^{1/2}(1+\pi(\gamma^2-1))^{3/2}}{2(\pi(1-\pi))^{3/2}} xN^{1/2}.
\end{equation}

This allows us then to express
%
%
\begin{eqnarray}
&&P\bigl(h_t(yt)\geq\bigl((1-2\eta)y+2\eta(1-\eta)\bigr)t - \sqrt{(1-\pi
-\eta)}t^{1/2}x\bigr) \nonumber\\[-8pt]\\[-8pt]
&&\qquad= P\bigl(L_2(N,M)\leq t\bigr),\nonumber
\end{eqnarray}
where $t$ is as above. It then follows after a little algebra that
Theorem \ref{LPP_fluctuations} applies and gives that these
probabilities asymptotically equal
%
%
\begin{equation}
F_G \biggl(\frac{x}{2\sqrt{\pi(1-\pi)}} \biggr)F_G
\biggl(\frac{x}{2\sqrt{\eta(1-\eta)}} \biggr),
\end{equation}
as desired to prove this part of Theorem \ref{thm_conj}.

\textit{$F_0$ case}: similarly to the previous case, we rehash the
height function event in terms of the LPP with two-sided boundary
conditions event and show that the desired asymptotic probabilities
arise from Theorem \ref{LPP_fluctuations}. The region of
$\rho_-,\rho_+$ for which we wish to prove $F_0$ fluctuations
corresponds to $\pi+\eta>1$. We wish to prove $F_0$ fluctuations for
all $y$ such
that $1-2\rho_-<y<1-2\rho_+$. Translating this region into
$\eta,\pi,\gamma$ variables exactly corresponds to the region in
which Theorem~\ref{LPP_fluctuations} implies $F_0$ fluctuations.

We wish to compute the asymptotic formula for
%
%
\begin{equation}
P \bigl(h_t(yt)\geq\tfrac{1}{2}(y^2+1)t -
2^{-1/3}(1-y^2)^{2/3}t^{1/3}x \bigr),
\end{equation}
where we have used the fact that $\bar{h}(y) = (y^2+1)/2$ in the region
of $y,\rho_-,\rho_+$ which we are considering.
Without loss of generality, let us assume that $y\geq0$ (the other
case follows similarly). As before, set
%
%
\begin{eqnarray}
N-M &=& yt,
\\
N+M &=& \tfrac{1}{2}(y^2+1)t - 2^{-1/3}(1-y^2)^{2/3}t^{1/3}x.
\end{eqnarray}

The height function event we are considering has the same
probability as
%
%
\begin{equation}
P\bigl(L_2(N,M)\leq t\bigr).
\end{equation}

Using the equations for $N-M$ and $N+M$, we can solve for
%
%
\begin{eqnarray}
N&=&\tfrac{1}{4}(1+y)^2t- 2^{-4/3}(1-y^2)^{2/3}t^{1/3}x,
\\
M&=&\tfrac{1}{4}(1-y)^2t- 2^{-4/3}(1-y^2)^{2/3}t^{1/3}x.
\end{eqnarray}

These expressions allow us to express $\gamma$ asymptotically as
%
%
\begin{eqnarray}\quad
\gamma^2 &=& \biggl(\frac{1-y}{1+y} \biggr)^2
\biggl[1+2^{2/3}xt^{-2/3}(1-y^2)^{2/3} \biggl(\frac{1}{(1+y)^2}-\frac
{1}{(1-y)^2} \biggr) \biggr]\nonumber\\[-8pt]\\[-8pt]
&&{}+o(t^{-2/3}).\nonumber
\end{eqnarray}

We may use Lemma \ref{inverting} to invert our expression
for $M$ in terms of $t$. Asymptotically
%
%
\begin{equation}\label{t_for_M}
t=\frac{4M}{(1-y)^2} + 2^{4/3}\frac
{(1+y)^{2/3}}{(1-y)^{2}}xM^{1/3}+o(M^{1/3}).
\end{equation}

This implies that $t^{-2/3} = \frac{4^{-2/3}}{(1-y)^{-4/3}}M^{-2/3}
+o(M^{-2/3})$. This can be plugged into our expression for $\gamma^2$
and gives a new expression for $\gamma^2$ in term of $M$ now:
%
%
\begin{equation}
\gamma^2 = \biggl(\frac{1-y}{1+y} \biggr)^2 \biggl[1-2^{4/3}xM^{-2/3}\frac
{y}{(1+y)^{4/3}} \biggr] +o(M^{-2/3}).
\end{equation}

From this, we may find that
%
%
\begin{equation}
y=\frac{1-\gamma}{1+\gamma} -\frac{(1-\gamma)\gamma}{(1+\gamma
)^{5/3}}xM^{-2/3} +o(M^{-2/3}).
\end{equation}

This can then be substituted into (\ref{t_for_M}) which gives
%
%
\begin{equation}t=(1+\gamma^{-1})^2 M + \frac{(1+\gamma
)^{4/3}}{\gamma
}M^{1/3}x + o(M^{1/3}).
\end{equation}

Plugging this into $P(L_2(N,M)\leq t)$, we find that our height
function probability is asymptotically equal to
%
%
\begin{equation}
P \biggl(L_2(N,M)\leq(1+\gamma^{-1})M + \frac{(1+\gamma
)^{4/3}}{\gamma}M^{1/3}x +o(M^{1/3}) \biggr).
\end{equation}
As it was already noted, the $y,\rho_-,\rho_+$ which we are considered
maps exactly onto the range of $\eta,\pi,\gamma$ which the LPP
probability above is asymptotically equal to $F_0(x)$ and hence the
same holds for the the height function probability.
\end{pf*}

\section*{Acknowledgments}
This paper came out of discussions at the Courant ASEP seminar. The
authors would like to thank Percy Deift for jointly organizing this
seminar with them, as well as thank all of the participants. The authors
appreciate helpful discussion of this material with Antonio
Auffinger. The authors are grateful to Tomohiro Sasamoto for catching
an important calculation error in an early reading of this paper, and
to Timo Sepp\"{a}l\"{a}inen for orienting them to what has been
previously proved.

%

%
\printaddresses

\end{document}